\documentclass[10pt]{amsart}

\usepackage{preamble}

\begin{document}

\title{Ext groups in Homotopy Type Theory}



\author{J. Daniel Christensen}
\address{University of Western Ontario, London, Ontario, Canada}
\email{jdc@uwo.ca}

\author{Jarl G.\ Taxerås Flaten}
\address{University of Western Ontario, London, Ontario, Canada}
\email{jtaxers@uwo.ca}

\date{August 15, 2025}

\maketitle

\begin{abstract}
  Ext groups are fundamental homological invariants which have important applications
  in homotopy theory and algebra.
  In particular, they appear in the classical universal coefficient theorem,
  a key computational tool in homotopy theory.
  Motivated by the goal of extending such tools to synethetic homotopy theory,
  we develop the theory of Yoneda Ext groups~\cite{Yoneda1954} over a ring in homotopy type theory (HoTT) and describe their interpretation into an \(\infty\)-topos.
  The Yoneda approach to Ext groups does not require projective or injective resolutions,
  which is a crucial in HoTT since we do not know that such resolutions exist.
  While it produces group objects that are a priori \emph{large},
  we show that the $\Ext^1$ groups are equivalent to small groups,
  leaving open the question of whether the higher Ext groups are essentially small as well.
  We also show that the $\Ext^1$ groups take on the usual form as a product of cyclic groups whenever the input modules are finitely presented and the ring is a PID (in the constructive sense).

  When interpreted into an \(\infty\)-topos of sheaves on a 1-category, our Ext groups recover (and give a resolution-free approach to) \emph{sheaf} Ext groups, which arise in algebraic geometry~\cite{Gro57}. (These are also called ``local'' Ext groups.)
  We may therefore interpret results about Ext from HoTT and apply them to sheaf Ext.
  To show this, we prove that injectivity of modules in HoTT interprets to internal injectivity in these models.
  It follows, for example, that sheaf Ext can be computed using resolutions which are projective or injective in the sense of HoTT, when they exist, and we give an example of this in the projective case.
  We also discuss the relation between internal \(\Zb G\)-modules (for a \(0\)-truncated group object \(G\)) and abelian groups in the slice over \(BG\), and study the interpretation of our Ext groups in both settings.
\end{abstract}

\vspace{-2em}
\tableofcontents

\section{Introduction}
\label{sec:intro}

We begin the study of homological algebra in homotopy type theory (HoTT) by
developing the theory of Ext groups of modules over a ring $R$.
Ext groups are important algebraic invariants, and also have
many applications in homotopy theory.
Classically, Ext groups are ingredients in the universal coefficient theorem
for cohomology, and we hope to use the results here to obtain a
universal coefficient spectral sequence in homotopy type theory.
In addition, the results discussed here were used in~\cite{BCFR}
to show that certain types must be products of Eilenberg--Mac Lane spaces.

It is common to define Ext groups in terms of injective or projective resolutions, as these often exist in classical settings.
However, we do not know whether such resolutions generally exist in HoTT.
(For instance, the usual way of constructing projective resolutions uses that free modules are projective, which is not provable in HoTT.)
We therefore follow a less common approach to defining Ext groups which avoids resolutions, namely the approach of Yoneda~\cite{Yoneda1954, Yoneda1960}.
This approach is also described in~\cite{Mac63}.
As we will see, carrying out this approach in HoTT is not straightforward, and univalence will play a key role in our definition of \(\Ext^1\).

In this approach, given modules $A$ and $B$ over a ring $R$, the \(n\)-th Ext group \(\eExt{R}^n(B,A)\) is defined as the set of path components of the space of length-\(n\) exact sequences
\[ 0 \to A \to E_1 \to E_2 \to \cdots \to E_n \to B \to 0 \]
of $R$-modules.
For \(n = 1\), this definition can be elegantly carried out in HoTT, by virtue of univalence.
Specifically, we define \(\SES_R(B,A)\) to be the type of short exact sequences from \(A\) to \(B\) (\cref{dfn:ses1}).
Using univalence, we see that paths in \(\SES_R(B,A)\) correspond to isomorphisms between short exact sequences,
so our type is capturing the correct notion.
We define (\cref{dfn:ext-1})
\[ \Ext_R^1(B,A) \ \defeq \ \Tr{\SES_R(B,A)}_0 . \]
The definition of \(\eExt{R}^n\) for \(n > 1\) is more difficult in HoTT, because we do not know how to correctly represent the space of length-\(n\) exact sequences.
Instead we define $\Ext_R^n(B,A)$ to be the set-quotient of a certain type \(\ES{R}^n\) equipped with a relation (\cref{dfn:ext-n}).
This approach is described in~\cite{Mac63}, and has been formalized in HoTT in~\cite{Fla23} for \(R = \Zb\).
We show that these types are abelian groups through an operation known as the \emph{Baer sum}~\cite{Baer}.

One aspect of these resolution-free definitions of Ext groups is that they produce
types lying in a larger universe.
We show that, for $n=1$, our Ext groups are essentially small:

\theoremstyle{plain}
\newtheorem*{thm:ses-baut}{\cref{thm:ses-baut}}
\begin{thm:ses-baut}
  Let \(B\) and \( A \) be abelian groups. We have a natural equivalence
    \[ \SES_\Zb (B,A) \ \simeq \ \big( \K{B}{2} \pto \K{A}{3} \big) .  \]
  In particular, \( \SES_\Zb (B,A) \) and \( \Ext_\Zb^1(B,A) \) are equivalent to small types.
\end{thm:ses-baut}

Here \(\pto\) denotes the type of pointed maps.
This result easily implies that \(\SES_R(B,A)\) is essentially small for modules over a general ring \(R\) (\cref{cor:ext1-small}).
However, we don't know whether the higher Ext groups are essentially small in general.
From this theorem, we also deduce the usual six-term exact sequences of Ext groups from a fibre sequence, in the special case \(R \jeq \Zb\) (\cref{pro:cov-six-term-ab,pro:contra-six-term-ab}).

The usual long exact sequences also exist, and the contravariant one has been formalized for $\Zb$ in~\cite{Fla23}.
Using these, we show in \cref{pro:ext-projective-resolution} that our Ext groups can be computed using projective (and injective) resolutions, whenever one is at hand.
It follows that our Ext groups yield right-derived functors of the hom-functor of modules whenever one has enough projectives or injectives.
More generally, we show:

\newtheorem*{thm:ext-delta-universal}{\cref{thm:ext-delta-universal}}
\begin{thm:ext-delta-universal}
  The large \(\delta\)-functor \(\lbrace \Ext_R^n(-,A) \rbrace_{n : \Nb}\) is universal, for any \(R\)-module \(A\).
\end{thm:ext-delta-universal}

The definition of a universal \(\delta\)-functor is recalled in \cref{dfn:delta-functor}.

We stress that the higher Ext groups \(\Ext_\Zb^n(B,A)\) need not vanish for \(n > 1\) in our setting.
Indeed, there are models of HoTT in which these are non-trivial, as we discuss below.
Nevertheless, we show that these Ext groups do vanish whenever \(R\) is a (constructive) PID and the module \(B\) is finitely presented (\cref{cor:ext-vanishes-fp}).
Moreover, when \(A\) is also finitely presented, then we get the usual description of \(\Ext_R^1(B,A)\) as a product of cyclic groups (\cref{pro:ext-fp}).

For a (0-truncated) group \(G\), it is well-known that HoTT lets one work \(G\)-equivariantly by working in the context of the classifying type \(BG\).
In an \(\infty\)-topos, this corresponds to working in the slice over the object \(BG\).
We study our constructions from this perspective in \cref{ssec:Ext-ZG}, which later lets us work out concrete examples of the interpretation of our Ext groups in \cref{ssec:Ext-BG}.
An abelian group ``in the context of \(BG\)'' is simply a map \(BG \to \Ab\).
Since the type of modules is 1-truncated, it is equivalent to replace \(BG\) with a pointed, connected type \(X\) (and \(G\) by \(\pi_1(X)\)).
To emphasize that our proofs do not require any truncation assumptions, we choose to work with such an \(X\).
We show that the category \(X \to \Ab\) is equivalent to the category of \emph{\(\Zb \pi_1(X)\)-modules} (\cref{pro:zg-mod}), where \(\Zb \pi_1(X)\) is the usual group ring.
When working in the context of \(X\), we carry out operations pointwise.
For example, given \(B, A : X \to \Ab\), we form the ``\(\Omega X\)-equivariant'' type of short exact sequences
\( x \mapsto \SES_\Zb(B_x, A_x) : X \to \Type. \)
Of course, we can also consider the type \(\SES_{\Zb \pi_1(X)}(B_\pt,A_\pt)\) of short exact sequences of \(\Zb \pi_1(X)\)-modules,
where \(B_\pt\) and \(A_\pt\) are the \(\Zb \pi_1(X)\)-modules corresponding to the families \(B\) and \(A\).
These are related:

\newtheorem*{thm:Z-vs-ZG}{\cref{thm:Z-vs-ZG}}
\begin{thm:Z-vs-ZG}
  For any \(B, A : X \to \Ab\), we have an equivalence
  \[ \prod_{x : X} \SES_{\Zb}(B_x,A_x) \ \simeq \ \SES_{\Zb \pi_1(X)}(B_\pt,A_\pt). \]
\end{thm:Z-vs-ZG}

We deduce the usual formula relating \(\eExt{}^1\) and cohomology with local coefficients:

\newtheorem*{cor:Ext-local-cohomology}{\cref{cor:Ext-local-cohomology}}
\begin{cor:Ext-local-cohomology}
  For any \(M : X \to \Ab\), we have a group isomorphism
  \[ H^1(X; M) \ \simeq \ \Ext_{\Zb \pi_1(X)}^1(\Zb, M_\pt), \]
  where the left-hand side is the cohomology of \(X\) with local coefficients in \(M\), and \(\Zb\) on the right has trivial \(\Zb \pi_1(X)\)-action.
\end{cor:Ext-local-cohomology}

In \cref{sec:interpretation}, we interpret our main results and constructions from HoTT into an \(\infty\)-topos \(\Xc\).
Given a ring \(R\) in an \(\infty\)-topos, it was shown in~\cite[Theorem~4.3.4]{Fla22} that the interpretation of the category of (``small'') \(R\)-modules from HoTT yields an internal category in \(\Xc\) which represents the presheaf sending an object \(X \in \Xc\) to the category of \(({X{\times}R})\)-modules in the slice \(\Xc / X\).
(Here \((X{\times}R)\) is a ring object in this slice.)
Building on this, in \cref{ssec:SES} we show that the object of short exact sequences \( \iSES{R}(B,A) \) between two modules \( A \) and \( B \) in \( \Xc \) represents the presheaf
\[ X \longmapsto \eSES{(X{\times}R)}(X \times B, X \times A) \ : \ \Xc^\op \longrightarrow \spaces , \]
where \(\eSES{(X{\times}R)}\) denotes the (1-truncated) space of short exact sequences between (``small'') \((X{\times}R)\)-modules.
From this description we show how to recover the classical Ext groups in \cref{cor:recover-Ext}.

An interesting result is that in certain \(\infty\)-toposes, the interpretation of our Ext groups recovers \emph{sheaf} Ext (\cref{dfn:sheaf-Ext}), which has been studied in algebraic geometry~\cite{Gro57}.
A consequence of this is that we can study sheaf Ext via the internal logic of a (higher) topos, and indeed the statements we prove for Ext in HoTT can be interpreted to give results for sheaf Ext.
Moreover, this shows that our definitions generalize Yoneda's approach from ordinary Ext groups to sheaf Ext.

The precise theorem is:

\newtheorem*{thm:recover-sheaf-ext}{\cref{thm:recover-sheaf-ext}}
\begin{thm:recover-sheaf-ext}
  Suppose sets cover in \( \Xc \).
  For any \( X \in \Xc \), ring \( R \in \Xc / X \) and \( R \)-module \( B \), the functor \( \iExt{R}^n(B,-) : \eMod{R} \to \eAb{\Xc / X} \) is naturally isomorphic to the sheaf Ext functor $\sExt{R}^n(B,-)$.
\end{thm:recover-sheaf-ext}

The meaning of ``sets cover'' is that any object admits an effective epimorphism from a \(0\)-truncated object (\cref{dfn:n-types-cover}).
Sets cover in any \(\infty\)-topos of \(\infty\)-sheaves on a 1-category.

Sheaf Ext is traditionally defined using injective resolutions (which always exist in these models), however our definition does not rely on the existence of enough injectives.
To prove this theorem we show that injectivity in HoTT interprets to \emph{internal injectivity} in these \(\infty\)-toposes (\cref{cor:hott-internal-injectivity}), which in turn follows from showing that internal injectivity is stable by base change in these models.
Our proof of stability uses (and partly generalizes) results of Roswitha Harting~\cite[Theorem~1.1]{Har83} (for abelian groups) and Blechschmidt~\cite[Proposition~3.7]{Ble18} (for modules) which show that internal injectivity of modules is stable by base change in any elementary 1-topos.
In addition, \cite[Theorem~3.8]{Ble18} shows that (externally) injective modules are always internally injective, which means that our Ext groups can be computed using the same resolutions used for sheaf Ext.

We also study various notions of projectivity of modules in \(\Xc\), namely the usual \emph{(external) projectivity}, \emph{internal projectivity}, and the notion of projectivity from HoTT.
In order to understand the relation between these notions, we provide examples which demonstrate that neither of external and internal projectivity imply the other (\cref{exa:external-not-internal-projective,exa:internal-not-external-projective}).
The example of an internally projective module that is not externally projective is an adaptation of an argument by Todd Trimble.
Moreover, we show that free modules on internally projective objects satisfy the notion of projectivity of modules from HoTT (\cref{pro:hott-projective-free}).
Using this fact, we demonstrate that our higher Ext groups need not vanish even over \(\Zb\) by computing a nontrivial \(\Ext_\Zb^2\).
There are also known computations of sheaf Ext which demonstrate this.

Finally, in \cref{ssec:Ext-BG} we study the theory developed throughout \cref{sec:interpretation} in some concrete situations.
In particular, we relate our Ext groups of abelian groups in a slice \(\Xc / BG\) to Ext groups of abelian groups in the base \(\Xc\) (\cref{pro:zg-ext-slice}), and deduce a vanishing result (\cref{cor:zg-ext-vanish}).
We also generalize another result of Harting~(\cref{pro:internally-injective-G-module}), and discuss the connection between our Ext groups in the slice \(\spaces / BG\) and ordinary Ext groups of \(\Zb G\)-modules (\cref{exa:ext-G-modules}).

\addtocontents{toc}{\SkipTocEntry}
\subsection*{Formalization}
Many of the results from \cref{sec:ext} have been formalized in HoTT and some have been contributed to the Coq-HoTT library~\cite{coqhott} under the namespace {\tt HoTT.Algebra.AbSES}.
This includes the (contravariant) six-term exact sequence and the Baer sum, the latter of which was formalized with the help of Jacob Ender.

Results about the higher Ext groups, including their construction, is currently in the separate repository \href{https://github.com/jarlg/Yoneda-Ext}{github.com/jarlg/Yoneda-Ext}.
The main result contained therein is the (contravariant) long exact sequence of Ext groups associated to any short exact sequence.
More information can be found in~\cite{Fla23}.

All of our formalization has been done for Ext groups over the ring \(\Zb\) for pragmatic reasons.

\addtocontents{toc}{\SkipTocEntry}
\subsection*{Open questions}
We list some outstanding questions.
\begin{enumerate}
\item In HoTT, is the abelian group \(\Ext_R^n(B,A)\) equivalent to a small type for \(n \geq 2\)?
  Is it independent of the universe for \(n \geq 2\)?
  (The case \(n = 1\) is answered by \cref{thm:ses-baut}.)
\item In HoTT, are injectivity and projectivity of $R$-modules (\cref{dfn:HoTT-inj,dfn:HoTT-proj}) independent of the universe?
\item In an \(\infty\)-topos, do HoTT-injectivity and HoTT-projectivity of $R$-modules only depend on the 1-topos of 0-truncated objects?
Do they agree with internal injectivity and internal projectivity?%
\footnote{David Wärn has now shown that internal projectivity does \emph{not} imply HoTT-projectivity; see \cite[Theorem~7]{Warn24}.}
These would follow from proving that internal injectivity and internal projectivity are pullback stable.
\item Does the interpretation of \(\Ext_R^n(B,A)\) into an \(\infty\)-topos depend only on the 1-topos of 0-truncated objects?
  (For \(\infty\)-toposes in which sets cover, this is answered by \cref{thm:recover-sheaf-ext}.)
\end{enumerate}

\addtocontents{toc}{\SkipTocEntry}
\subsection*{Notation and conventions}
Our setting is Martin-Löf type theory with a hierarchy of univalent universes, as in the HoTT Book~\cite{hottbook}, whose notation we generally follow.
Univalence and function extensionality are taken as axioms, and are frequently used implicitly.
The only HITs we require are pushouts and truncations.
All of our groups, rings and modules are assumed to be sets (i.e., $0$-truncated).
We write \(\Type\) for a fixed universe, and \(\Type_*\) for the universe of pointed types.
\cref{sec:interpretation} has its own paragraph on notation.

\addtocontents{toc}{\SkipTocEntry}
\subsection*{Acknowledgements}
We thank Jacob Ender for contributions to the formalization of the Baer sum,
which are described in~\cite{Fla23} and have been contributed to~\cite{coqhott}.
We thank David Wärn for a question which led to \cref{ssec:fp}.
This article has been improved thanks to very helpful feedback from referees.
We acknowledge the support of the Natural Sciences and Engineering Research Council of Canada (NSERC), RGPIN-2022-04739.

\section{Ext in HoTT}
\label{sec:ext}

In this section, we develop the theory of Yoneda Ext groups in HoTT.
Many of the results we show have classical analogues, in which case our contribution is the verification that these results hold in our setting as well.
Nevertheless, our proofs and definitions make use of univalence and truncations, and we make constructive considerations (particularly in \cref{ssec:fp}), all of which do not feature in the traditional theory.

Let \(R\) be a ($0$-truncated) ring throughout this entire section.

\subsection{The type of short exact sequences}
Fix two left \( R \)-modules \( A \) and \( B \) throughout this section.
Below we define the type \( \SES_R(B,A) \) whose elements are short exact sequences
\[ 0 \to A \xra{i} E \xra{p} B \to 0 \]
in \( \Mod{R} \).
The type \( \SES_R(B,A) \) is a $1$-type, and we define the set \( \Ext_R^1(B,A) \) of extensions to be its set-truncation.
By characterizing paths in \( \SES_R(B,A) \), we will show that an extension is trivial if and only if it is \emph{merely split}.

A homomorphism of \(R\)-modules is an \define{epimorphism} (resp.\ a \define{monomorphism}) if and only if its underlying function is surjective (resp.\ an embedding).
We write \(\Epi{R}(E,B)\) and \(\Mono{R}(A,E)\) for the set of \(R\)-module epimorphisms and monomorphisms, respectively.

\begin{dfn}
  Let \( A \xra{i} E \xra{p} B \) be two composable homomorphisms in \( \Mod{R} \).
  Whenever the composite \( p \circ i \) is trivial, there is a unique induced homomorphism \( i' : A \to \mathsf{ker}(p) \).
  If \(i'\) is an epimorphism, then \(i\) and \(p\) are \textbf{exact}:
  \[ \mathsf{IsExact}(i,p) \defeq  \sig{h : \prod_{a:A} p(i(a)) = 0}{\mathsf{IsEpi}(i') }. \]
\end{dfn}

\begin{dfn} \label{dfn:ses1}
The \define{type of short exact sequences from \(A\) to \(B\)} is:
  \[ \SES_R(B,A) \defeq \sum_{E : \Mod{R}}\ \sum_{i : \Mono{R}(A,E) }\ \sum_{p : \Epi{R}(E,B)} \mathsf{IsExact}(i,p) . \]
  We often denote a short exact sequence by its middle module $E$,
  and write \(i_E : A \hra E \) and \(p_E : E \thra B\) for the inclusion and projection homomorphisms.
  The type is pointed by the split short exact sequence \(A \xra{\mathsf{in}_A} A \oplus B \xra{\mathsf{pr}_B} B \).
\end{dfn}

As defined, \( \SES_R \) quantifies over \( \Mod{R} \) and is therefore a \emph{large} type.
It is moreover a \(1\)-type, since \( \Mod{R} \) is a \(1\)-type,
$i$ and $p$ range over sets, and $\mathsf{IsExact}(i,p)$ is a proposition.
This mirrors the classical fact that the category of module extensions of $B$ by $A$---whose maps are homomorphisms \( E \to E' \) making the relevant triangles commute---is a groupoid.
The following proposition strengthens this connection:

\begin{pro} \label{pro:path-ses1}
  For short exact sequences \(A \xra{i} E \xra{p} B\) and \(A \xra{j} F \xra{q} B\), we have
  \[ (E =_{\SES_R(B,A)} F) \ \simeq \ \sig{\phi : \Mod{R}(E,F)}{(\phi \circ i = j) \land (p = q \circ \phi)} . \]
\end{pro}

\begin{proof}
  It follows from the characterization of paths in \(\Sigma\)-types and transport in function types that
  \[ (E =_{\SES_R(B,A)} F) \ \simeq \ \sum_{\phi : E \cong F} {(\phi \circ i = j) \land (p = q \circ \phi)} , \]
  where \(E \cong F\) denotes \(R\)-module isomorphisms.
  The stated equivalence now follows from the short five lemma.
\end{proof}

This lets us compute the loop space of \( \SES_R(B,A) \) as in~\cite{Retakh1986}.

\begin{cor} \label{cor:loops-ses1}
  We have a natural isomorphism \( \widehat{(-)} : \Omega \SES_R(B,A) \simeq \Mod{R}(B,A) \) of groups.
\end{cor}

\begin{proof}
  By the previous proposition, a path \( A \oplus B =_{\SES_R} A \oplus B \) corresponds to a homomorphism \( \phi : A \oplus B \to A \oplus B \) which respects \( \mathsf{in}_A \) and \( \mathsf{pr}_B \).
  Thus we get an \(R\)-module homomorphism \( \widehat{\phi} : B \to A \) as the composite
  \[ \begin{tikzcd}[cramped]
      \widehat{\phi} : B \rar["\mathsf{in}_B"] & A \oplus B \rar["\phi"] & A \oplus B \rar["\mathsf{pr}_A"] & A.
    \end{tikzcd} \]

  Conversely, for any homomorphism \( f : B \to A \), we can define the homomorphism
  \[ (a,b) \mapsto (a + f(b), b) : A \oplus B \to A \oplus B \]
  which respects \( \mathsf{in}_A \) and \( \mathsf{pr}_B \).
  These associations are easily shown to be mutually inverse, yielding a bijection \( \Omega \SES_R(B,A) \simeq \Mod{R}(B,A) \).
  To see that it's an isomorphism of groups, consider a composite path \( \phi \cdot \psi \).
  The associated \(R\)-module homomorphism \( A \oplus B \to A \oplus B \) is given by the composite
  \[ (a, b) \longmapsto (a + \widehat{\phi}(b), b) \longmapsto (a + \widehat{\phi}(b) + \widehat{\psi}(b), b). \]
  Hence \( \widehat{\phi \cdot \psi} = \widehat{\phi} + \widehat{\psi} \), as required.
\end{proof}

In \cite{Mac63}, Mac Lane produces the set (underlying the abelian group) of extensions by applying \( \pi_0 \) to the groupoid of short exact sequences.
We now do the corresponding thing:

\begin{dfn} \label{dfn:ext-1}
  The \define{set of extensions of $B$ by $A$} is
  \( \Ext_R^1(B,A) \defeq \Tr{\SES_R(B,A)}_0 \).
\end{dfn}

The following proposition characterizes trivial extensions.

\begin{pro} \label{pro:ext-trivial-iff-merely-splits}
  Let \( E \) be a short exact sequence from \(A\) to \( B \).
  Then \( E \) is trivial in \( \Ext_R^1(B,A) \) if and only if \( p \) \emph{merely splits}, i.e., the following proposition holds:
  \[ \Bigtr{ \sig{s : \Mod{R}(B,E) } {p \circ s = \id_B } } . \]
\end{pro}

\begin{proof}
  First of all, by the characterization of paths in truncations~\cite[Theorem~7.3.12]{hottbook} we have
  \[ \big( \tr{E}_0 =_{\Ext_R^1(B,A)} 0 \big) \ \simeq \ \Tr{ E =_{\SES_R(B,A)} A \oplus B } . \]
Forgetting about truncations, the right-hand side holds if and only if \( p \) splits, by the usual argument. This in turn implies the statement on the truncations.
\end{proof}

We conclude this section by showing that \(\Ext_R^1\) defines a bifunctor which lands in abelian groups.
This is also shown in~\cite[Section~3.2]{Fla23}, but we give a different proof.

\begin{dfn} \label{dfn:pullback-pushout-ses}
  Let \( A \to E \to B \) be a short exact sequence of \( R \)-modules.
  \begin{enumerate}[(i)]
  \item For \( f : A \to A' \), the \define{pushout \(f_*(E)\) of \(E\) along \( f \)} is the short exact sequence defined by the dashed homomorphisms below:
    \[ \begin{tikzcd}
        A \rar \dar["f"] \ar[dr, phantom, "\ulcorner" at end] & E \rar \dar & B \period \\
        A' \rar[dashed] & f_*(E) \ar[ur, dashed, bend right]
      \end{tikzcd} \]
    Here the curved arrow is defined to make the triangle commute and to be zero on $A'$.
  \item For \( g : B' \to B \), the \define{pullback \( g^*(E) \) of \( E \) along \( g \)} is the short exact sequence defined by the dashed homomorphisms below:
    \[ \begin{tikzcd}
          & g^*(E) \dar \rar[dashed] \ar[dr, phantom, "\lrcorner" at start] & B' \dar["g"] \\
        A \rar \ar[ur, dashed, bend left] & E \rar & B \period
      \end{tikzcd} \]
    Here the curved arrow is defined to make the triangle commute and to be zero into $B'$.
  \end{enumerate}
  For any \(R\)-module \(M\), these operations define functions
  \[ f_* : \SES_R(M,A) \to \SES_R(M,A') \quad \text{and} \quad g^* : \SES_R(B,M) \to \SES_R(B',M) . \]
\end{dfn}

The pushout and pullback operations commute (in the sense that \(f_*g^*(E) = g^*f_*(E)\) whenever this expression makes sense) meaning \(\Ext_R^1(-,-)\) is a bifunctor into \(\Set\).
This has been formalized in~\cite[Proposition~9]{Fla23}.
Before making this bifunctor land in \(\Ab\), we also need to \emph{detect} pushouts (and pullbacks) of short exact sequences, in the following sense.

\begin{lem} \label{lem:ses-pushout}
  Suppose given a diagram
  \[ \begin{tikzcd}
      A \rar \dar["\alpha"] & E \rar \dar & B \dar["\beta"] \\
      A' \rar & F \rar & B
  \end{tikzcd} \]
  with short exact rows.
  If \(\beta = \id_B\), then there is a path \(\alpha_*(E) = F\) of short exact sequences.
  \qed
\end{lem}

The dual statement for pullbacks requires the leftmost vertical homomorphism to be the identity.
For a proof, the reader may consult~\cite[Proposition~7]{Fla23} or Lemmas~III.1.2 and III.1.4 in~\cite{Mac63}.

That \(\Ext_R^1(B,A)\) is an abelian group follows from the following proposition.

\begin{pro} \label{pro:ext-additive}
  The contravariant functor \(\Ext_R^1(-,A)\) takes arbitrary (set-indexed) coproducts to products, and the covariant functor \(\Ext_R^1(B,-)\) preserves finite products.
\end{pro}

\begin{proof}
  We first show that \(\Ext_R^1(-,A)\) takes arbitrary coproducts to products.
  To that end, let \(X\) be a set and consider a family \(B : X \to \Mod{R}\).
  Theorem~3.3.10 of~\cite{Fla22} produces an exact coproduct functor \(\bigoplus_X : (X \to \Mod{R}) \to \Mod{R}\)
  from the category of $X$-indexed families of $R$-modules to $R$-modules.
  We will show that the natural map \(\phi : \Ext^1_R(\bigoplus_{x : X} B_x, A) \to \prod_{x : X} \Ext^1_R(B_x, A)\) is a bijection for any \(R\)-module \(A\).
  This natural map is defined as follows.
  Since we are defining maps between sets, we may pick representatives of extensions.
  Given a short exact sequence
  \[  A \to E \to \bigoplus_{x : X} B_x ,\]
  define \(E_x\) to be the result of pulling back \(E\) along the natural map \(B_x \to \bigoplus_X B\) for \(x : X\).
  A map in the inverse direction is given as follows.
  A family
  \( (A \to F_x \to B_x)_{x : X} \)
  of short exact sequences yields a short exact sequence
  \[ \bigoplus_X A \to \bigoplus_{x:X} F_x \to \bigoplus_{x:X} B_x \]
  by exactness of \(\bigoplus_X\), which by pushing out along \(\nabla : \bigoplus_X A \to A\) yields an element of \(\Ext^1_R(\bigoplus_X B, A)\).

  Starting from a short exact sequence \(A \to E \to \bigoplus_X B\), the following diagram exhibits the bottom row as the pushout of the top row by \cref{lem:ses-pushout}, showing that \(\phi\) is a section:
  \[ \begin{tikzcd}
      \bigoplus_X A \rar \dar["\nabla"] & \bigoplus_{x : X} E_x \dar \rar & \bigoplus_{x : X} B_x \dar[equals] \\
      A \rar & E \rar & \bigoplus_{x : X} B_x \period
    \end{tikzcd} \]
  Here the middle vertical arrow is induced from the maps \( (E_x \to E)_{x : X}\) coming from the definition of \(E_x\) as a pullback.

  Similarly, for any family \( F \defeq (A \to F_x \to B_x)_{x : X}\), the following diagram exhibits the top row as the pullback of the bottom row along \(B_y \to \bigoplus_X B\), for any \(y : X\), since the composite of the left vertical maps is the identity on $A$:
  \[ \begin{tikzcd}
      A \dar \rar & F_y \dar \rar & B_y \dar \\
      \bigoplus_X A \dar["\nabla"] \rar & \bigoplus_{x : X} F_x \rar \dar & \bigoplus_{x : X} B_x \dar[equals] \\
      A \rar & \nabla_* \bigoplus_{x : X} F_x \rar & \bigoplus_{x : X} B_x \period
    \end{tikzcd} \]
  Here the maps from the top row to the middle row are given by the inclusion of the \(y\)-summand.
  This shows that \(\phi\) is a retraction, hence a bijection.

  To show that \(\Ext^1_R(B,-)\) preserves finite products, it suffices to check that it preserves the empty product and binary products.
  The former is clear, so we proceed to handle binary products.
  The natural map
  \( \Ext^1_R(B, A_0 \oplus A_1) \to \Ext^1_R(B,A_0) \times \Ext^1_R(B,A_1) \)
  is given by pushing out along the two projections of \(A_0 \oplus A_1\).
  To get a map in the opposite direction, we take the biproduct of the given extensions
  (using that biproducts are exact)
  and then pull back along \(\Delta : B \to B \oplus B\).
  Showing that these two maps are inverses
  is straightforward.
\end{proof}

\begin{cor}\label{cor:baer-sum}
  Let \(A\) and \(B\) be \(R\)-modules.
  The set \(\Ext_R^1(B,A)\) is naturally an abelian group.
\end{cor}

\begin{proof}
  We have just shown that the functor \(\Ext_R^1(B,-) : \Mod{R} \to \Set\) preserves finite products.
  It follows that it preserves group objects.
  But any \(R\)-module is itself an abelian group object in \(\Mod{R}\) (in a unique way), so we are done.
\end{proof}

The binary operation on \(\Ext_R^1(B,A)\) is called the \emph{Baer sum}.
A concrete description of this operation, which has been formalized in joint work with Jacob Ender, is discussed in~\cite[Section~3.3]{Fla23}.
We also mention that if the ring \(R\) is commutative, then \(\Ext_R^1(B,A)\) is naturally an \(R\)-module.

We also record the following lemma for later use.
\begin{lem} \label{lem:pullback-pushout-iso}
  Let \(f : A \to A'\) and \(g : B' \to B\) be isomorphisms of \(R\)-modules.
  For any short exact sequence \(A \xra{i} E \xra{p} B\), we have
  \( g^* f_*(E) = (A' \xra{i \circ f^{-1}} E \xra{g^{-1} \circ p} B'). \) \qed
\end{lem}

\subsection{Classifying extensions and smallness of \texorpdfstring{\( \Ext^1 \)}{Ext¹}}
We remarked after Definition~\ref{dfn:ses1} that \( \SES_R \) is a large type, and consequently \( \Ext^1_R \) is a large abelian group.
This is not surprising, since our definition mirrors that of the external Yoneda Ext groups in an abelian category, and examples of abelian categories are known where these are proper classes.
However, our \( \Ext_R^1 \) groups turn out to be equivalent to small types.

In this and the next sections we will be working with fibre sequences, whose definition we now recall.

\begin{dfn}
  A \define{fibre sequence} consists of a sequence \( F \xpra{i} X \xpra{f} Y \) of pointed maps equipped with a pointed homotopy witnessing that \( f \circ i \) is constant such that the induced map \( F \to \fib{f} \) is an equivalence.
\end{dfn}

Recall that~\cite[Theorem~5.1]{BvDR} produces the following equivalence of categories for \( n \geq 2 \):
\[ \K{-}{n} : \Ab \ \simeq \ \{ \t{ pointed, $(n{-}1)$-connected $n$-types } \} : \Omega^n . \]
A (pointed) map \(f\) between \((n{-}1)\)-connected \(n\)-types is \((n{-}1)\)-connected if and only if \(\Omega^n(f)\) is \((-1)\)-connected (i.e., an epimorphism), so these classes of maps correspond under this equivalence.
Moreover, like any equivalence of categories, this equivalence preserves pullback squares.
By combining these facts, we see that the following classes of squares correspond under this equivalence:

\[ \begin{tikzcd}
  A \rar \dar \ar[dr, phantom, "\lrcorner" at start] & E \dar[two heads] && \K{A}{n} \rar \dar \ar[dr, phantom, "\lrcorner" at start] & \K{E}{n} \dar \\
  0 \rar & B && \pt \rar & \K{B}{n}
  \end{tikzcd} \]
Here, a square of the kind on the left is precisely a short exact sequence \( A \to E \to B \),
and a square of the kind on the right precisely says that $K(A,n) \to K(E,n) \to K(B,n)$ is
a fibre sequence.
In other words, short exact sequences are sent to fibre sequences, and vice-versa, under the stated equivalence.
We use this to prove the following:

\begin{thm}\label{thm:ses-baut}
  Let \(B\) and \( A \) be abelian groups. We have a natural equivalence
    \[ \SES_\Zb (B,A) \ \simeq \ \big( \K{B}{2} \pto \K{A}{3} \big) .  \]
  In particular, \( \SES_\Zb (B,A) \) and \( \Ext_\Zb^1(B,A) \) are equivalent to small types.
\end{thm}

The right-hand side of the equivalence above is moreover equivalent to \(\big( \K{B}{n} \pto \K{A}{n+1} \big) \) for \(n \geq 2\), since \(\Omega\) is an equivalence in this range.
(See~\cite[Theorem~6.7]{BvDR}, with their $n = 0$ and their $k$ equal to our $n$.
Their $F$ is our $\Omega$.)
The right-hand side is also equivalent to \(\big( \K{B}{n} \pto \BAut{\K{A}{n}} \big)\),
since $\K{A}{n+1}$ is the 1-connected cover of $\BAut{\K{A}{n}}$
(\cite{shulman_emspace} and~\cite[Proposition 5.9]{BCFR}).

\begin{proof}
  We define maps in both directions and show that they are mutual inverses.
  To go from left to right, we apply \(\K{-}{3}\) to a short exact sequence \(A \to E \to B\) to get a fibre sequence, then we negate the maps and take the fibre:
  \[ \begin{tikzcd}[column sep=large]
      \K{B}{2} \rar[dashed] & \K{A}{3} \rar["-\K{i}{3}"]  & \K{E}{3} \rar["-\K{p}{3}"] & \K{B}{3} \period
    \end{tikzcd} \]
  The fibre is naturally equivalent to \(\K{B}{2}\), as displayed, since the three rightmost terms form a fibre sequence.
  This process yields a map from left to right.

Conversely, given a map \(f : \K{B}{2} \pto \K{A}{3}\), we get a fibre sequence
\[ \K{A}{2} \longrightarrow F \longrightarrow \K{B}{2} \xra{\ f\ } \K{A}{3}, \]
where $F$ is a pointed, 1-connected 2-type.
Taking loop spaces twice yields a short exact sequence \(A \to \Omega^2 F \to B\)
of abelian groups.

We first consider the composite starting and ending at \(\SES_\Zb(B,A)\).
Starting with a short exact sequence \(A \to E \to B\), we apply \(\K{-}{3}\),
negate the maps and take three fibres, producing the sequence
\[ \begin{tikzcd}[column sep=large]
    \K{A}{2} \rar["\K{i}{2}"] & \K{E}{2} \rar["\K{p}{2}"] & \K{B}{2} \period
  \end{tikzcd} \]
Here we have used that taking three fibres negates maps, by~\cite[Lemma~8.4.4]{hottbook}.
We then apply $\Omega^2$, which yields the original short exact sequence, since \(\Omega^2 \circ \K{-}{2}\) is the identity.

For the composite starting and ending at \(\big(\K{B}{2} \pto \K{A}{3}\big)\), we will use that the map
\[ \Omega : \big(\K{B}{3} \pto \K{A}{4}\big) \longrightarrow \big(\K{B}{2} \pto \K{A}{3}\big) \]
is an equivalence.
Write \(B\) for the inverse.
This equivalence implies that any map \(\phi : \K{B}{2} \pto \K{A}{3}\) fits into the following fibre sequence:
\[ \begin{tikzcd}
    \K{B}{2} \rar["\phi"] & \K{A}{3} \rar & \fib{-B\phi} \rar & \K{B}{3} \rar["-B\phi"] & \K{A}{4} \period
  \end{tikzcd} \]
Applying \(\Omega^3\) to the middle three terms produces the short exact sequence associated to \(\phi\), but with the maps negated.
Since \(\Omega^3\) is an equivalence and commutes with negation of maps, this means that the middle three terms are equal to \(\K{-}{3}\) applied to the short exact sequence associated to \(\phi\) with the maps negated.
It immediately follows that the composite starting and ending at \(\big(\K{B}{2} \pto \K{A}{3}\big)\) is equal to the identity as well, so the maps we defined are mutual inverses.

It is straightforward to check that the map from right to left is natural in both \(B\) and \(A\).
\end{proof}

\cref{thm:ses-baut} is about \emph{abelian} extensions of abelian groups.
Results similar to \cref{thm:ses-baut}, but for central extensions,
appear in~\cite{BvDR, djm, NSS14, Sco20}, where one takes $n = 1$.

\begin{cor} \label{cor:ext1-small}
  Let \( B \) and \( A \) be $R$-modules.
  Then both \( \SES_R (B,A) \) and \( \Ext^1_R (B,A) \) are equivalent to small types.
\end{cor}

\begin{proof}
  Let \( U : \Mod{R} \to \Ab \) denote the forgetful functor.
  The fibre $\fib{u}(E)$ of the induced map \( u : \SES_R (B,A) \to \SES_\Zb (UB, UA) \) over an extension \( E \) of abelian groups is equivalent to the type of
  $R$-module structures on $E$ which make $i$ and $p$ into $R$-module homomorphisms.
  Since the type of $R$-module structures on $E$ and the two conditions are all small,
  this fibre is equivalent to a small type.
  Now we use that \( \SES_R (B,A) \simeq \sum_{E : \SES_\Zb (UB, UA)} \fib{u}(E) \), where the latter is equivalent to a small type.
\end{proof}

\begin{rmk}\label{rmk:not-a-product1}
Classically, one argument that the external Ext group $\eExt{R}^1(B, A)$ is small is that the underlying
set of any extension $E$ of $A$ by $B$ is isomorphic to the product set $A \times B$.
However, this can fail in models of HoTT, such as in the Sierpi\'nski $\infty$-topos,
as we show in \cref{rmk:not-a-product}.
If we allow more general situations, smallness of external Ext can also fail.
For example,~\cite{wofsey} describes a locally small abelian category in which the external
Yoneda Ext group $\eExt{\Zb}^1(Z, Z)$ can be a proper class for a certain object $Z$.
We believe that this category can arise as the category of abelian group
objects in an elementary $\infty$-topos of $G$-spaces, where $G$ is the free abelian group
on a proper class of generators, and for each object, all but a set of generators
are required to act trivially.
In this setting, $Z$ is the interpretation of the integers, and so
the interpretation of our $\Ext_{\Zb}^1(Z, Z)$ is zero, since the integers are projective in the sense of HoTT (see \cref{dfn:HoTT-proj,pro:HoTT-proj}).
This illustrates that it is somewhat surprising that the interpretation of $\Ext^1_R(B,A)$ is small
in every model of HoTT.
\end{rmk}

\begin{rmk}
It follows from the equivalence $\SES_\Zb(B,A) \simeq ( \K{B}{2} \pto \K{A}{3} )$
in \cref{thm:ses-baut} that $\SES_{\Zb}(B, A)$ is independent of the choice of
universe containing $A$ and $B$.
Therefore, the same holds for $\Ext^1_{\Zb}(B, A)$.
The argument in the proof of \cref{cor:ext1-small} shows that these statements
are also true when $\Zb$ is replaced by a general ring $R$,
since the set of $R$-module structures on an abelian group $E$ is independent
of the choice of universe.
\end{rmk}

\subsection{The six-term exact sequences}
\label{ssec:six-term-Z}

For $A \to E \to B$ a short exact sequence of $R$-modules and $M$ another $R$-module,
there are covariant and contravariant six-term exact sequences of abelian groups
\[
  0 \to \Mod{R}(M, A) \xra{i_*} \Mod{R}(M, E) \xra{p_*} \Mod{R}(M,B)
    \to \Ext_R^1(M,A) \xra{i_*} \Ext_R^1(M,E)  \xra{p_*} \Ext_R^1(M,B)
\]
and
\[
             \Ext_R^1(A,M) \xla{i^*} \Ext_R^1(E, M) \xla{p^*} \Ext_R^1(B,M)
  \leftarrow \Mod{R}(A, M) \xla{i^*} \Mod{R}(E, M) \xla{p^*} \Mod{R}(B,M)
  \leftarrow 0 .
\]
These can be proved following Theorem~3.2 of \cite{Mac63}, and the contravariant version has been formalized in~\cite[Proposition~19]{Fla23}.
Here we find it interesting to give different arguments in the special
case where $R \jeq \Zb$.

\begin{pro} \label{pro:cov-six-term-ab}
  Let \( A \xra{i} E \xra{p} B : \Ab \) be a short exact sequence, and \( M : \Ab \).
  Then
  \[ \SES_\Zb (M,A) \xra{i_*} \SES_\Zb (M,E) \xra{p_*} \SES_\Zb (M,B) \]
  is a fibre sequence, where the maps are given by pushing out.
\end{pro}

\begin{proof}
  Applying \( \K{-}{3} \) to the given short exact sequence produces a fibre sequence
  \[ \K{A}{3} \to \K{E}{3} \to \K{B}{3} . \]
  Since $(Z \pto -)$ preserves fibre sequences for any pointed type $Z$, we can apply $(\K{M}{2} \pto -)$
  to obtain a fibre sequence
  \[ \big(\K{M}{2} \pto \K{A}{3}\big) \to \big(\K{M}{2} \pto \K{E}{3}\big) \to \big(\K{M}{2} \pto \K{B}{3}\big), \]
  where the maps are given by post-composition.
  \cref{thm:ses-baut} then gives the desired fibre sequence, since naturality means that post-composition corresponds to pushout of short exact sequences.
\end{proof}

Using \cref{cor:loops-ses1}, one can show that the long exact sequence of homotopy groups associated to the fibre sequence above recovers the usual covariant six-term exact sequence of Ext groups mentioned at the beginning of this section.

We now give the dual result, which can similarly be shown to produce the contravariant six-term exact sequence of Ext groups.
The construction of the following fibre sequence is more difficult, because we need to map out of a fibre sequence (not into, as in the previous proposition).

\begin{pro} \label{pro:contra-six-term-ab}
  Let \( A \xra{i} E \xra{p} B : \Ab \) be a short exact sequence, and \( M : \Ab \).
  Then 
  \[ \SES_\Zb(A,M) \xla{i^*} \SES_\Zb(E,M) \xla{p^*} \SES_\Zb(B,M) \]
  is a fibre sequence, where the maps are given by pulling back.
\end{pro}

\begin{proof}
  Applying \( \K{-}{2} \) to the given short exact sequence produces a fibre sequence
  \[ \K{A}{2} \to \K{E}{2} \to \K{B}{2} . \]
  Let $C$ be the cofibre of the map $\K{A}{2} \to \K{E}{2}$, which comes
  with a natural map $C \to \K{B}{2}$.
  Since $(- \pto Z)$ sends cofibre sequences to fibre sequences for any $Z$,
  we can apply $(- \pto \K{M}{3})$ to obtain a fibre sequence
  \[ \big(\K{A}{2} \pto \K{M}{3}\big) \leftarrow \big(\K{E}{2} \pto \K{M}{3}\big) \leftarrow \big(C \pto \K{M}{3}\big) . \]
  We claim that $\big(C \pto \K{M}{3}\big) \simeq \big(\K{B}{2} \pto \K{M}{3}\big)$,
  from which the statement follows, as in the proof of the previous result.
  Since $\K{M}{3}$ is a 3-type, it suffices to prove that $\Tr{C}_3 \; \simeq \Tr{\K{B}{2}}_3$,
  and for this it suffices to show that the map $C \to \K{B}{2}$ is 3-connected,
  using~\cite[Lemma~7.5.14]{hottbook}.
  The map $C \to \K{B}{2}$ is the cogap map associated to the map $\K{E}{2} \to \K{B}{2}$
  and the base point inclusion $1 \to \K{B}{2}$.
  Since $\K{B}{2}$ is connected, it suffices to check the connectivity of the fibre of this map
  over the base point.
  By~\cite[Theorem~2.2]{Rij17}, this fibre is the join $\K{A}{2} * \Omega \K{B}{2}$
  of the fibres, which is (1+0+2)-connected, as required.
  (This fact about connectivities of joins is proved in~\cite[\texttt{Join.v}]{coqhott}.
  It also follows from~\cite[Corollary~2.32]{CS20},
  since the join is the suspension of the smash product.)
\end{proof}

\subsection{Higher Ext groups}
\label{ssec:les}

The definition of higher Ext groups from~\cite[Chapter~XII]{Mac63} or~\cite[pp.~216]{Yoneda1954} can be translated to HoTT and has already been formalized (for \(R \jeq \Zb\), but the arguments work for a general ring) in~\cite{Fla23} along with the contravariant long exact sequence.
An account of the covariant long exact sequence that can be carried out in our setting may be found in~\cite[Chapter~VII.5]{Mit65}.
We first discuss the definition of \(\Ext_R^n\), referring the reader to~\cite{Fla23} for further details.
The long exact sequence of Ext groups \cite[Theorem~26]{Fla23} makes the collection \(\big\lbrace \Ext_R^n(-,A) \big\rbrace_{n : \Nb}\) into a \emph{(large) \(\delta\)-functor}, for any \(A\).
In \cref{thm:ext-delta-universal}, we show that this \(\delta\)-functor is universal, as expected.

\begin{dfn} \label{dfn:es-n}
  Let \(B\) and \(A\) be \(R\)-modules.
  The type \(\ES{R}^n(B,A)\) is inductively defined to be
  \[ \ES{R}^n(B,A) \defeq
    \begin{cases}
      \; \Mod{R}(B,A) &\text{if } n \jeq 0, \\[4pt]
      \; \SES_R (B,A) &\text{if } n \jeq 1, \\[4pt]
      \displaystyle
      \sum_{C : \Mod{R}} \ES{R}^m(C,A) \times \SES_R(B,C) & \text{if } n \jeq m+1, m > 0.
    \end{cases} \]
  There is an evident \textbf{splicing} operation
  \( \splice : \ES{R}^m(C,A) \times \SES_R(B,C) \to \ES{R}^{m+1}(B,A) \) for any \(C : \Mod{R}\), which is given by pushing out along a homomorphism when \(m \jeq 0\).
\end{dfn}

Our splicing operation is written in diagrammatic order, as in~\cite{Mac63}.
An element of \(\ES{R}^n\) consists of \(n\) short exact sequences which can be spliced in succession from left to right.
It is straightforward to define a more general splicing operation where the right factor can have arbitrary length.
The base point of \(\ES{R}^0(B,A)\) is the zero homomorphism,
the base point of \(\ES{R}^1(B,A)\) is the split short exact sequence, and for $n \geq 2$
the base point of \(\ES{R}^n(B,A)\) is recursively defined to be
\[ 0 \to A \to A \to 0 \to \cdots \to 0 \to B \to B \to 0 , \]
where there are \(n-2\) intermediate occurrences of \(0\)
(and just a zero map when $n = 2$).

The type $\ES{R}^n(B,A)$ doesn't correctly represent the type of length $n$ exact
sequences.  For example, one cannot define $\Ext_R^n(B,A)$ as the 0-truncation of this space.
Even more, we would like to have that $\ES{R}^{n+1}(B,A)$ is a delooping of $\ES{R}^n(B,A)$
for each $n$, as in~\cite{Retakh1986}, but this does not hold.
However, we can define
$\Ext_R^n(B,A)$ as the set-quotient of $\ES{R}^n(B,A)$ by a certain relation.

\begin{dfn}\label{dfn:esrel}
  Let \(n : \Nb\).
  For \(E, F : \ES{R}^n(B,A) \) define a (type-valued) relation inductively by
  \[ E \esrel F \defeq
    \begin{cases}
      \hspace{1em} E = F &\text{if } n \jeq 0, 1 \\[4pt]
      \displaystyle
      \sum_{\beta : \Mod{R}(C, C')} (E_l \esrel \beta^*(F_l)) \times (\beta_*(E_r) = F_r)
            &\text{if } n > 1, E \jeq (C, E_l, E_r) \text{ and } F \jeq (C', F_l, F_r).
    \end{cases} \]
\end{dfn}

This relation is clearly reflexive, and one easily shows by induction that it is transitive.
We now define the higher Ext groups as the set-quotient of \(\ES{R}^n\) by this relation.

\begin{dfn} \label{dfn:ext-n}
  Let \(B\) and \(A\) be \(R\)-modules.
  For \(n : \Nb\), define the \textbf{set of length-\boldmath{\(n\)} extensions of \boldmath{\(B\)} by \boldmath{\(A\)}} to be
  \[ \Ext_R^n(B,A) \defeq
    \begin{cases}
      \Mod{R}(B,A) &\text{if } n \jeq 0, \\[4pt]
      \Ext^1_R(B,A) &\text{if } n \jeq 1, \\[4pt]
      \Tr{\ES{R}^n(B,A) / {\esrel}}_0 &\text{if } n > 1.
    \end{cases} \]
\end{dfn}

The splicing operation respects the relation \(\esrel\) and thus passes to the quotient \(\Ext_R^n\).
The same is true for pushouts and pullbacks of length-\(n\) exact sequences, which makes \(\Ext_R^n\) into a profunctor.

In~\cite[III.5]{Mac63}, Mac Lane considers ``ladders'' of homomorphisms between length-\(n\) exact sequences, and defines two such sequences to be \emph{congruent} if there exists a zig-zag of such ladders beween them.
He defines \(\Ext_R^n\) by quotienting out by this congruence.
To connect our definition of \(\Ext_R^n\) to Mac Lane's, and for later use, we state the following definition and lemma.
For an element \(E : \ES{R}^n(B,A)\), we write \(A_i \xra{i_{E_i}} E_i \xra{p_{E_i}} A_{i+1} \) for the \(i\)-th splice factor, with \( 1 \leq i \leq n \), \(A_1 \jeq A\) and \(A_{n+1} \jeq B\).

\begin{dfn}
  Let \(E, E' : \ES{R}^n(B,A)\) for \(n \geq 1\).
  A \define{morphism \(f : E \to E'\)} consists of a family \(f_i : E_i \to E'_i \) of $R$-module homomorphisms for \( 1 \leq i \leq n\) such that \(f_1 \circ i_{E_1} = i_{E'_1}\), \(p_{E_n} = p_{E'_n} \circ f_n \), and the following outer square commutes for \(1 \leq i < n\):
  \[ \begin{tikzcd}[sep=large]
      E_i \dar["f_i"] \rar["p_{E_{i}}"] & A_{i+1} \dar[dashed, "\hat{f}_{i+1}"] \rar["i_{E_{i+1}}"] & E_{i+1} \dar["f_{i+1}"] \\
      E'_i \rar["p_{E'_{i}}"] & A'_{i+1} \rar["i_{E'_{i+1}}"] & E'_{i+1} .
      \end{tikzcd} \]
    We also define \(\hat{f}_{i+1}\) to be the displayed restriction of \(f_{i+1}\), which makes the inner squares commute.
\end{dfn}

\begin{lem} \label{lem:ladder-homo-ext-path}
  Let \( n \geq 1\), and let \(E, E' : \ES{R}^n(B,A)\).
  The types \(E \esrel E'\) and \(E \to E'\) are logically equivalent.
\end{lem}

\begin{proof}
  The base case \(n \jeq 1\) is \cref{pro:path-ses1}, so we need only show that the statement holds for \(n+1\) supposing it holds for \(n \geq 1\).
  Write \(E \jeq E_l \splice E_r\) and \(E' \jeq E_l' \splice E_r'\).

  We first construct a map \((E \esrel E') \to (E \to E')\) as follows.
  Given \((\beta, \rho, q)\) in the domain, define \(f_{n+1}\) to be the composite of the natural morphism \( E_r \to \beta_*(E_r) \) and the morphism \(\beta_*(E_r) \to E'_r\) underlying \(q\).
  By recursing on \(\rho\) we get a morphism \(E_l \to \beta^*(E'_l)\), which we compose with the natural morphism \(\beta^*(E'_l) \to E'_l\) to get \(f_{1, \ldots, n} : E_l \to E'_l\).
  The collection \((f_i)_{1 \leq i \leq n + 1}\) then defines a morphism \(E \to E'\), by construction.

  For the other direction, suppose we have a morphism \(f : E \to E'\).
  This gives a morphism $\hat{f}_{n+1} : A_{n+1} \to A'_{n+1}$, using the notation
  from the previous definition.
  This morphism will play the role of $\beta$ from \cref{dfn:esrel}.
  By \cref{lem:ses-pushout} we have that \((\hat{f}_{n+1})_* E_r = E'_r\), and we can use the induction hypothesis to get an element of
  \( E_l \esrel  \hat{f}^*_{n+1} E'_l \), completing the proof.
\end{proof}

\begin{rmk}
  The previous lemma implies that the set-quotients of \(\ES{R}^n(B,A)\) given by the relations `\(\esrel\)' and `\(\to\)' are equivalent.
  For the formalisation we found it easier to work with the relation `\(\esrel\)'.

  More abstractly, one can inductively define \(\Ext_R^{n+1}(B,A)\) as the set-coequalizer of
  \[ \begin{tikzcd}
     \displaystyle
      \sum_{C, C' : \Mod{R}} \, \sum_{f : C \to C'} \Ext_R^n(B, C) \times \Ext^1(C', A) \ar[rr, shift left=1] \ar[rr, shift right=1]
      && \displaystyle \sum_{C : \Mod{R}}  \Ext_R^n(B, C) \times \Ext_R^1(C,A) \comma
    \end{tikzcd} \]
  where one map is given by pushing out along \(f\) and the other is given by pulling back along \(f\). The reader may consult~\cite[p.~553]{Yoneda1960} for details about this definition.
\end{rmk}

We define a Baer sum on \(\Ext_R^{n+2}(B,A)\) by \(E + F \defeq \Delta^* \nabla_* (E \oplus F) \), and this makes \(\Ext_R^{n+2}(B,A)\) into an abelian group, for all \(n : \Nb\).
We now show that the splicing operation is a homomorphism:

\begin{pro} \label{pro:splice-homomorphism}
  Let \(F : \Ext_R^n(B,A)\) for some \(n \geq 0\), and let \(C\) be an \(R\)-module.
  The splicing operation \(F \splice - : \Ext_R^1(C, B) \to \Ext_R^{n+1}(C,A) \) is a homomorphism.
\end{pro}

\begin{proof}
  For \(n \jeq 0\) the statement is that pushing out an extension along the homomorphism \(F\) is a homomorphism, which follows from \cref{cor:baer-sum}.
  Suppose \(n \geq 1\), and let \(E, E' : \Ext_R^1(C,B)\).
  An easy consequence of \cref{lem:ladder-homo-ext-path} is that \(\nabla^* F = \nabla_* (F \oplus F) \), which we use to deduce:
  \begin{align*}
    F \splice\, (E + E') \
    &\jeq \ F \splice \Delta^* \nabla_* (E \oplus E') &\text{definition of the Baer sum} \\
    &\jeq \ \Delta^* \big( F \splice \nabla_* (E \oplus E') \big) &\text{since pullback acts on the right-hand factor}\\
    &= \ \Delta^* \big( \nabla^* F \splice\, (E \oplus E') \big) &\text{using the definition of \(\esrel\)}\\
    &= \ \Delta^* \big( \nabla_* (F \oplus F) \splice\, (E \oplus E') \big) &\text{using that \(\nabla^* F = \nabla_* (F \oplus F) \)} \\
    &\jeq \ \Delta^* \nabla_* \big( (F \oplus F) \splice\, (E \oplus E') \big) &\text{since pushout acts on the left-hand factor}\\
    &= \ \Delta^* \nabla_* \big( (F \splice E) \oplus (F \splice E') \big) &\text{since \(\splice\) and \(\oplus\) commute} \\
    &\jeq \ (F \splice E) + (F \splice E') . &&\qedhere
  \end{align*}
\end{proof}

A similar argument shows that splicing is a homomorphism in the other variable.
We use these facts to show the following:

\begin{cor}
For each $n$, $\Ext_R^n(-, -)$ is an additive functor in each variable.
\end{cor}

\begin{proof}
This is well-known for $n=0$ and was proved for $n=1$ in \cref{pro:ext-additive}.
We'll prove that $\Ext_R^{n+1}(-,A)$ is additive for $n \geq 1$; the other case is similar.
First we show that for $f, g : \Mod{R}(B', B)$,
$\Ext_R^{n+1}(f+g,A) = \Ext_R^{n+1}(f,A) + \Ext_R^{n+1}(g,A)$
as maps $\Ext_R^{n+1}(B, A) \to \Ext_R^{n+1}(B', A)$.
This is a proposition, so we can assume we have a length $n+1$ sequence of
the form $F \splice E$ for $F : \ES{R}^n(C, A)$ and $E : \SES_R(B, C)$.
We then have the following equalities in \(\Ext_R^{n+1}(B',A)\):
\begin{align*}
  (f+g)^* (F \splice E) \
  &\jeq \ F \splice \, (f+g)^* E &\text{since pullback acts on the right-hand factor} \\
  &= \ F \splice \, (f^* E + g^* E) &\text{since $\Ext_R^1$ is additive} \\
  &= \ (F \splice f^*E) + (F \splice g^* E) &\text{by \cref{pro:splice-homomorphism}} \\
  &\jeq \ f^*(F \splice E) + g^*(F \splice E) .
\end{align*}
This shows that
\( \Ext_R^{n+1}(f+g,A) = \Ext_R^{n+1}(f,A) + \Ext_R^{n+1}(g,A) \).
It is also straightforward to see that $\Ext_R^{n+1}(-, A)$ sends the zero
map to the zero map.
It follows from this that $\Ext_R^{n+1}(-,A)$ preserves direct sums, since a
direct sum $B_1 \oplus B_2$ is characterized by having homomorphisms
$i_k : B_k \to B_1 \oplus B_2$ and $p_k : B_1 \oplus B_2 \to B_k$ for $k = 1, 2$
such that $p_1 i_1 = \id$, $p_2 i_2 = \id$, $p_1 i_2 = 0$, $p_2 i_1 = 0$
and $i_1 p_1 + i_2 p_2 = \id$,
and these equations are preserved by $\Ext_R^{n+1}(-, A)$.
\end{proof}

Our next goal is to show that, for each $R$-module $A$, the collection \(\big\lbrace \Ext_R^n(-,A) \big\rbrace_{n : \Nb}\) has the structure of a \emph{(large) universal \(\delta\)-functor}, a concept that we define now.

\begin{dfn} \label{dfn:delta-functor}
  A \define{\(\delta\)-functor structure} on a collection of additive functors \(\lbrace T^n : \Mod{R}^\op \to \Ab \rbrace_{n : \Nb}\) associates to any short exact sequence \(A \to E \to B\) of \(R\)-modules a \define{connecting homomorphism} \(\delta^n_E : T^n(A) \to T^{n+1}(B) \) for each \(n : \Nb\), such that:
  \begin{enumerate}[(i)]
  \item  The following long complex is exact:
    \[ 0 \to T^0(B) \to T^0(E) \to T^0(A) \xra{\delta_E^0} T^1(B) \to \cdots \to T^n(E) \to T^n(A) \xra{\delta^n_E} T^{n+1}(B) \to \cdots . \]
  \item For any morphism of short exact sequences as on the left below, the square on the right commutes for every \(n : \Nb\):
    \[ \begin{tikzcd}
        A \rar \dar["\alpha"] & E \rar \dar & B \dar["\beta"] && T^n(A') \dar["T^n_\alpha"] \rar["\delta^n_F"] & T^{n+1}(B') \dar["T^{n+1}_\beta"] \\
        A' \rar & F \rar & B' && T^n(A) \rar["\delta^n_E"] & T^{n+1}(B) \period
      \end{tikzcd} \]
  \end{enumerate}
  A \define{\(\delta\)-functor}\footnote{More precisely, this is the notion of a \emph{contravariant}, \emph{cohomological} \(\delta\)-functor~\cite[Chapter~2.1]{Weibel}.} is such a collection equipped with a \(\delta\)-functor structure.
  Replacing \(\Ab\) above with the category \(\Ab'\) of large abelian groups%
  \footnote{For $A, B : \Ab'$, we still write $\Ab(A, B)$ for the large abelian group of group homomorphisms.},
  we obtain the notion of a \define{large \(\delta\)-functor}.

  If \(T\) and \(S\) are (large) \(\delta\)-functors, then a \define{morphism} \(f : T \to S\) of \(\delta\)-functors consists of a collection of natural transformations \(\lbrace f_n : T^n \Rightarrow S^n \rbrace_{n : \Nb}\) which respect the connecting homomorphisms.
  The (large) \(\delta\)-functor \(T\) is \define{universal} if the restriction map \((T \to S) \to (T^0 \Rightarrow S^0)\) is a bijection, for any (large) \(\delta\)-functor \(S\).
\end{dfn}

The splicing operation \(- \splice E \) defines connecting homomorphisms for the family
\(\lbrace \Ext_R^n(-,A) \rbrace_{n : \Nb} \) of contravariant functors,
and the long exact sequence from \cite[Theorem~26]{Fla23} shows that the first axiom holds.
It is straightforward to verify the second axiom.
Thus we have a large \(\delta\)-functor structure on \(\lbrace \Ext_R^n(-,A) \rbrace_{n : \Nb} \).
Below, we show that it is universal.
This fact is implicit in Yoneda's approach to \emph{satellites} in~\cite[Chapter~4]{Yoneda1960}, though he does not give an explicit proof of universality.
(Satellite is another word for (large) universal \(\delta\)-functor.)
However, Buchsbaum~\cite{Buchsbaum1960} constructs satellites which can be shown to be isomorphic to Yoneda's definition, and Buchsbaum does prove that his construction produces a universal \(\delta\)-functor (see his Proposition~4.3).

\begin{pro} \label{pro:ext-delta-map}
  Let \(T\) be a large \(\delta\)-functor and let \(A\) and \(B\) be \(R\)-modules.
  For each \(n : \Nb\), there is a homomorphism of abelian groups
  \( d_n : \Ext^n_R(B,A) \to \Ab\big( T^0(A), \, T^n(B) \big) \) which is natural in \(A\) and \(B\).
\end{pro}

\begin{proof}
  We proceed by induction on \(n\).
  Since \(T^0\) is an additive contravariant functor, it gives a homomorphism
  \[ \phi \longmapsto T^0_\phi \ : \ \Mod{R}(B,A) \longrightarrow \Ab \big(T^0(A),\, T^0(B)\big) \]
  which is natural in \(A\) and \(B\).
  We can therefore define \(d_0(\phi) \defeq T^0_\phi\).

  For \(n \jeq 1\), consider the map \(E \mapsto \delta_E^0 : \SES_R(B,A) \to \Ab\big(T^0(A),\, T^1(B)\big)\).
  Since the codomain is a set, we get our map \(d_1\) out of the set-truncation \(\Ext_R^1(B,A)\).
  We check that \(d_1\) is natural in \(A\); naturality in \(B\) is similar.
  Let \(f : A \to A'\) be a homomorphism.
  Our goal is to show that the square on the left commutes:
  \[ \begin{tikzcd}
      \Ext_R^1(B,A) \dar["f_*"] \rar["d_1"] & \dar["(T^0_f)^*"] \Ab\big(T^0(A),\, T^1(B)\big) && T^0(A') \dar["T^0_f"]  \rar["\delta^0_{f_*E}"] & T^1(B) \dar[equals] \\
      \Ext_R^1(B,A') \rar["d_1"] & \Ab\big(T^0(A'),\, T^1(B)\big) && T^0(A) \rar["\delta^0_E"] & T^1(B) \period
    \end{tikzcd} \]
  Since naturality is a proposition, we may pick an actual short exact sequence \(A \to E \to B\) in the top left corner.
  The question then is whether the equation \(d_1(f_*(E)) = d_1(E) \circ T^0_f\) holds.
  But this equation underlies the commuting square above on the right, which comes from part (ii) of the \(\delta\)-functor structure of \(T\) applied to the natural morphism \(E \to f_*E\) of short exact sequences.

  Now let $n \geq 1$ and assume that we have the natural homomorphism \(d_n\).
  We proceed to construct \(d_{n+1}\).
  First we define a map \(d'_{n+1} : \ES{R}^{n+1}(B,A) \to \Ab\big(T^0(A),\, T^n(B) \big) \) by
  \[ d'_{n+1}(F \splice E) \ \defeq \ \delta^n_E \circ d_n([F]) \ : \ T^0(A) \longrightarrow T^{n+1}(B), \]
  where \([F] : \Ext_R^n(C,A)\) is the equivalence class of \(F : \ES{R}^n(C,A)\).
  To descend \(d'_{n+1}\) to a map \(d_{n+1}\) on the quotient \(\Ext_R^{n+1}\), we need to show that it respects the relation on \(\ES{R}^{n+1}\).

  Suppose we have a relation \(E \esrel F\) in \(\ES{R}^{n+1}(B,A)\).
  Writing \(E \jeq (C, E_0, E_1)\) and \(F \jeq (C', F_0, F_1)\), the relation gives a map \(\beta : C \to C'\), a relation \(E_0 \esrel \beta^*(F_0)\) in \(\ES{R}^n(C,A)\), and a path \(\beta_*(E_1) = F_1\) in \(\SES_R(B',C)\).
  We need to argue that the outer square below commutes:
  \[ \begin{tikzcd}[sep=large]
      T^0(A) \rar["{d_n([F_0])}"] \dar["{d_n([E_0])}" swap] & \ar[dl, "T^n_\beta" swap] T^n(C') \dar["{\delta^n_{F_1}}"]  \\
      T^n(C) \rar["{\delta^n_{E_1}}"]& T^{n+1}(B) \period
    \end{tikzcd} \]
  The lower-right triangle commutes by condition (ii) of the \(\delta\)-structure of \(T\), using the map of short exact sequences \(E_1 \to F_1\) associated to the equality \(\beta_*(E_1) = F_1\).
  For the upper-left triangle, first note that we have \(d_n([E_0]) = d_n([\beta^*(F_0)])\).
  Naturality of \(d_n\) gives us further that \(d_n(\beta^*[F_0]) = T^n_\beta \circ d_n([F_0])\), from which we conclude that the upper-left triangle commutes.
  Thus we get the desired map \(d_{n+1}\) by passing to the quotient.

  It remains to show that \(d_{n+1}\) is natural and a homomorphism.
  By \cref{lem:naturality-implies-homomorphism} below, the latter follows from the former, so we only check naturality.
  First we check it in the first variable, so let \(f : A \to A'\) be a homomorphism.
  Since we need to show a proposition, we may consider an actual element \(F \splice E : \ES{R}^{n+1}(B,A)\).
  Then, since pushouts of longer exact sequences are defined recursively on the left factor, we have
  \begin{align*}
    d_{n+1}(f_* [F \splice E]) \
    &\jeq \ d_{n+1}([f_*(F) \splice E]) \
      \jeq \ \delta^{n+1}_E \circ d_n([f_*F]) \ \\
    &= \ \delta^{n+1}_E \circ d_n([F]) \circ T^0_f \ \jeq \ d_{n+1}([F \splice E]) \circ T^0_f ,
    \end{align*}
    where the only non-definitional equality uses naturality of \(d_n\).

    For naturality of \(d_{n+1}\) in the second variable, let \(g : B' \to B\) be a homomorphism.
    Again, we consider a general element \(F \splice E\) as above.
    Since pullback of longer exact sequences are defined directly on the last splice factor, we have
    \begin{align*}
      d_{n+1}(g^*[F \splice E]) \
      &\jeq \ d_{n+1}([F \splice g^*E]) \
        \jeq \ \delta^{n+1}_{g^*E} \circ d_n([F]) \ \\
      &= \ T^{n+1}_g \circ \delta^{n+1}_E \circ d_n([F])
      \ \jeq \ T^{n+1}_g \circ d_{n+1}([E \splice F]) ,
    \end{align*}
    where the only non-definitional equality comes from part (ii) of the \(\delta\)-functor structure of \(T\) applied to the natural morphism \(g^*E \to E\) of short exact sequences.
\end{proof}

The following is Proposition~4.1 in~\cite{Yoneda1960}, whose proof is easy to translate to our setting.
Abelian categories are defined as usual.

\begin{lem} \label{lem:naturality-implies-homomorphism}
  Let \(\cat{A}\) be an abelian category, and consider two additive functors \(S, T : \cat{A} \to \Ab\).
  Suppose \(\eta_A : S(A) \to T(A) \) is a collection of set-maps, natural in \(A \in \cat{A}\).
  Then each \(\eta_A\) is a homomorphism. \qed
\end{lem}

We come to the main result of this section.

\begin{thm} \label{thm:ext-delta-universal}
  The large \(\delta\)-functor \(\lbrace \Ext_R^n(-,A) \rbrace_{n : \Nb}\) is universal, for any \(R\)-module \(A\).
\end{thm}

\begin{proof}
  Let \(\lbrace T^n : \Mod{R}^\op \to \Ab' \rbrace_{n : \Nb}\) be a large \(\delta\)-functor.
  First note that \( \big(\Ext_R^0(-,A) \Rightarrow T^0\big) \simeq T^0(A) \), by the Yoneda lemma.
  (To be precise, we view \(\Mod{R}\) as an \(\Ab'\)-enriched category, and use the Yoneda lemma.)
  We will construct a morphism of (large) \(\delta\)-functors \( \Ext_R^n(-,A) \to T \) for any
  element \(\eta : T^0(A) \), and show that such morphisms are uniquely determined by their restriction to the zeroth level.

  Let \(n : \Nb\), let \(\eta : T^0(A)\), and let \(B\) be an \(R\)-module.
  Using \(d_n\) from the previous proposition, define
  \[ u_n(-) \jeq d_n(-, \eta) \ : \ \Ext^n_R(B,A) \longrightarrow T^n(B). \]
  Clearly \(u_n\) is a group homomorphism.
  Also, since \(u_0(\id_A) = T_{\id_A}(\eta) = \eta\), \(u_0\) corresponds to
  \(\eta\) under the Yoneda lemma.

  To see that \(\{u_n\}_{n : \Nb}\) is a morphism of \(\delta\)-functors, we need to show that it respects the connecting homomorphisms.
  To that end, let \(B \to E \to B'\) be a short exact sequence.
  We proceed by induction.
  For \(n \jeq 0\), we need to show that the following diagram commutes:
  \begin{equation} \label{eq:ext-delta-0}
    \begin{tikzcd}
      \Mod{R}(B,A) \rar["u_0"] \dar["- \splice E"] & T^0(B) \dar["\delta^0_E"]  \\
      \Ext^1_R(B',A) \rar["u_1"] & T^1(B') \period
    \end{tikzcd}
  \end{equation}
Let \(f : B \to A\) be an \(R\)-module morphism and
recall that $f \splice E \jeq f_*(E)$.
By definition, we have that \(u_1(f_*(E)) = \delta^0_{f_*(E)}(\eta)\).
Using functoriality of the \(\delta\)-structure of \(T\), the natural map of short exact sequences from \(E\) to \(f_*(E)\) yields a commuting square
\[ \begin{tikzcd}
    T^0(A) \dar["\delta^0_{f_*(E)}"] \rar["T_f"] & T^0(B) \dar["\delta^0_E"] \\
    T^1(B') \rar[equals] & T^1(B').
  \end{tikzcd} \]
Thus \(\delta^0_{f_*(E)}(\eta) = \delta^0_E(T_f(\eta))\).
The right-hand side is precisely \(\delta^0_E(u_0(f))\), concluding the base case.

For the inductive step, we need to show that the following square commutes,
for $n \geq 1$:
\begin{equation} \label{eq:ext-delta-n}
  \begin{tikzcd}
    \Ext^n_R(B,A) \rar["u_n"] \dar["- \splice E"] & T^n(B) \dar["\delta^n_E"] \\
    \Ext^{n+1}_R(B',A) \rar["u_{n+1}"] & T^{n+1}(B') \period
  \end{tikzcd}
\end{equation}
Whether the square commutes is a proposition, so we may choose a representative \(F\) of an element in the top left corner.
But then the square clearly commutes by the definition of \(u\) (and \(d\)).

It remains to show uniqueness of the \(\delta\)-functor morphism \(u\).
Specifically, we need to show that for any \(\delta\)-functor morphism \(\lbrace v_n : \Ext^n_R(-,A) \Rightarrow T^n \rbrace_{n : \Nb}\) such that \(v_0 = u_0\),
we have that \(v = u\).
To show that \(v_1(E) = u_1(E)\) for any \(E : \Ext^1_R(B,A)\), we may assume \(E\) is a short exact sequence.
Then we may consider diagram~\eqref{eq:ext-delta-0}, but with the bottom horizontal map being \(v_1\). 
For this \(E\), the top-left corner is \(\Mod{R}(A,A)\) and we may chase \(\id_A\) around the two sides of the square.
Since the square commutes, we get \(v_1(E) = \delta^0_E(\eta)\), and the right-hand side is \(u_1(E)\), by definition.
Similarly, for the inductive step we may write a general element of \(\Ext^{n+1}_R(B,A)\) as a splice \(F \splice E\) and consider diagram~\eqref{eq:ext-delta-n} with the lower horizontal map being \(v_{n+1}\).
(By the induction hypothesis the top horizontal map is \(u_n = v_n\).)
Chasing \(F\) around the two sides of the square, we get \(v_{n+1}(F \splice E) = \delta^n_E(u_n(E)) = u_{n+1}(E)\), as desired.
\end{proof}

\subsection{Computing Ext via projective resolutions}
\label{ssec:ext-proj-res}
In this section, we use the long exact sequence to show that our Ext groups can be computed using projective resolutions.
A dual argument shows the same for injective resolutions.
We begin by defining and characterizing projectivity and injectivity of modules in our setting.

\begin{dfn}\label{dfn:HoTT-proj}
We say that an $R$-module $P$ is \define{projective} if
for all $R$-modules $A$ and $B$ (in $\Type$),
every epimorphism $e : \Mod{R}(A,B)$
and every \( f : \Mod{R}(P,B) \), there merely exists a lift of $f$ through $e$:
    \[ \Bigtr{ \sig{ g : \Mod{R}(P,A) }{ e \circ g = f } }. \]
In other words, the post-composition map $e_* : \Mod{R}(P,A) \to \Mod{R}(P,B)$ is an epimorphism.
We write $\mathsf{IsProjective}(P)$ for this property.
\end{dfn}

It is clear that \(R^n\) is a projective \(R\)-module, for any ring \(R\) and natural number \(n\).
More generally, if $X$ is a projective \emph{set}, then the free $R$-module on $X$ is
a projective $R$-module.
In addition, binary coproducts of projective modules are easily seen to be projective.

The following reproduces a classical characterization of projective modules.

\begin{pro} \label{pro:HoTT-proj}
  Let \( P \) be an $R$-module. The following are equivalent:
  \begin{enumerate}[(i)]
  \item $P$ is projective.
  \item Every epimorphism \( p : \Mod{R}(A,P) \) merely splits, i.e., the following holds:
  \[ \Bigtr{ \sig{s : \Mod{R}(P,A)}{p \circ s = \id_P} }. \]
  \item \( \Ext_R^1(P, A) = 0 \) for all $R$-modules \( A \).
  \item \( \Ext_R^n(P, A) = 0 \) for all \(R\)-modules \(A\) and \(n \geq 1\).
  \end{enumerate}
\end{pro}

\begin{proof}
  The equivalence between $(i)$ and \( (ii) \) mirrors the classical argument, and
  the equivalence of $(ii)$ and $(iii)$ follows from Proposition~\ref{pro:ext-trivial-iff-merely-splits}.
  Clearly \((iv)\) implies \((iii)\); we will show the converse.

  Let \(n > 1\).
  Since we are proving a proposition, we may choose representatives and write a general element of \(\Ext_R^n(B,A)\) as
  \(\tr{(C, F, E)}_0\).
  By \((iii)\), we may assume that \(E\) is trivial.
  A relation \(0 \esrel F \splice E\) is then given by the constant map \(0 \to C\) and the fact that pushing out and pulling back along a constant map is trivial (as is easy to show by induction).
  Further details can be found in the formalization of~\cite{Fla23}.
\end{proof}

\begin{dfn}\label{dfn:HoTT-inj}
We say that an $R$-module $I$ is \define{injective} if
for all $R$-modules $A$ and $B$ (in $\Type$),
every monomorphism $m : \Mod{R}(A,B)$
and every \( f : \Mod{R}(A,I) \), there merely exists an extension of $f$ along $m$.
In other words, the pre-composition map $m^* : \Mod{R}(B,I) \to \Mod{R}(A,I)$ is an epimorphism.
We write $\mathsf{IsInjective}(I)$ for this property.
\end{dfn}

An argument dual to that given in \cref{pro:HoTT-proj} characterizes the injectives
using mere splittings or the condition that \( \Ext_R^1(B,I) = 0 \) for all \( B \).

In \cref{ssec:projectives}, we interpret projectivity into a model of HoTT and study its relation to existing notions of projectivity.
We do the same for injectivity in \cref{ssec:sheaf-Ext}.

Now we turn to computing \(\Ext_R^n\) from a projective resolution.
The argument is standard homological algebra, and the content is that it holds with the results available to us in homotopy type theory.
In the following, assume we have a projective resolution \( P_\bullet \) of \( B \).
This is equipped with a surjection \( p_0 : P_0 \to B \) inducing an isomorphism \( P_0 / \on{im}(P_1) = B \), and so \( P_1 \) surjects onto \( B_1 \defeq \on{ker}(p_0) \).
Continuing inductively, we may factor the projective resolution as follows: 
\begin{equation}\label{eq:resolution}
\begin{tikzcd}
{} \arrow[r, dotted]  & P_2 \arrow[d, two heads, "p_2"] \arrow[r] & P_1 \arrow[d, two heads, "p_1"] \arrow[r] & P_0 \arrow[d, two heads, "p_0"] \\
{} \arrow[ru, dotted] & B_2 \arrow[ru, hook, "i_2" swap]               & B_1 \arrow[ru, hook, "i_1" swap]               & B_0  \comma
\end{tikzcd}
\end{equation}
where \( B_0 \defeq B \) and \( B_{i+1} \defeq \on{ker}(p_i) \). Let \( P_{-1} \defeq 0 \) and \( i_0 \defeq 0 \) in the following.

\begin{pro} \label{pro:ext-projective-resolution}
  The abelian group \( \Ext_R^n(B,A) \) is the \( {n^{th}} \)cohomology of the cochain complex
  \[ \Mod{R}(P_\bullet, A) \defeq ( \cdots \to \Mod{R}(P_{n-1}, A)\to \Mod{R}(P_n, A) \to \Mod{R}(P_{n+1}, A) \to \cdots  ) . \]
\end{pro}

\begin{proof}
Applying $\Mod{R}(-, A)$ to \eqref{eq:resolution} gives a diagram
\[ \begin{tikzcd}
\ar[r] \ar[dr, "i_n^*" swap] \Mod{R}(P_{n-1}, A) &  \ar[r] \ar[dr, "{i_{n+1}^*}" swap] \Mod{R}(P_n, A)  & \Mod{R}(P_{n+1}, A) \\
     & \Mod{R}(B_n, A) \ar[u, hook, "p_n^*"] & \Mod{R}(B_{n+1}, A) \ar[u, hook, "{p_{n+1}^*}" swap]
  \end{tikzcd} \]
which has the chain complex across the top.
Since $p_{n+1}$ is an epimorphism, $p_{n+1}^*$ is a monomorphism, and so we get that
\( \on{ker}(\Mod{R}(P_n, A) \to \Mod{R}(P_{n+1}, A) ) = \on{ker}(i_{n+1}^*) \).
Since \( (B_{n+1} \xra{i_{n+1}} P_n \xra{p_n} B_n) \) is a short exact sequence,
the contravariant long exact sequence implies that \( \on{ker}(i_{n+1}^*) \) is
\( \Mod{R}(B_n, A) \), with \( p_n^* \) being the kernel inclusion.
Consequently,
\[ H^n (\Mod{R}(P_\bullet, A)) = \on{ker}(i_{n+1}^*) / \on{im}(\Mod{R}(P_{n-1},A)) = \Mod{R}(B_n, A) / \on{im}(i_n^*) = \on{cok}(i_n^*) . \]
If \(n = 0\), then this is \(\Mod{R}(B_0, A) \jeq \Ext_R^0(B,A) \), since \( i_0 \jeq 0 \).
To understand $\on{cok}(i_{n+1}^*)$, we use the full long exact sequence
\[ \begin{tikzcd}[row sep=small]
    0 \rar & \Mod{R}(B_n, A) \arrow[r, "p_n^*"]
    & \Mod{R}(P_n, A) \arrow[r, "i_{n+1}^*"] \ar[d, phantom, ""{coordinate, name=Z1}]
    & \Mod{R}(B_{n+1}, A) \ar[dll, rounded corners,
    to path={ -- ([xshift=2ex]\tikztostart.east)
      |- (Z1) [near end]\tikztonodes
      -| ([xshift=-2ex]\tikztotarget.west)
      -- (\tikztotarget)}] \\
    & {\Ext_R^1(B_n,A)} \arrow[r]
    & 0 \arrow[r] \ar[d, phantom, ""{coordinate, name=Z2}]
    & {\Ext_R^1(B_{n+1},A)} \ar[dll, rounded corners,
    to path={ -- ([xshift=2ex]\tikztostart.east)
      |- (Z2) [near end]\tikztonodes
      -| ([xshift=-2ex]\tikztotarget.west)
      -- (\tikztotarget)}] \\
    & {\Ext_R^2(B_n, A)} \arrow[r]
    & 0 \arrow[r] &
    {\Ext_R^2(B_{n+1},A)} \arrow[r,]    & \cdots
\end{tikzcd} \]
where the zeros down the middle column appear due to \( P_n \) being projective, by \cref{pro:HoTT-proj}.
It follows from the first connecting homomorphism that \( \on{cok}(i_{n+1}^*) = \Ext_R^1(B_n, A) \).
The subsequent connecting homomorphisms imply that \( \Ext_R^k(B_{n+1}, A) = \Ext_R^{k+1}(B_n, A) \) for \( k \geq 1 \).
Applying the second equality recursively gives \( \Ext_R^1(B_n, A) = \Ext_R^{n+1}(B,A) \),
for $n \geq 0$,
and so we conclude that $\on{cok}(i_{n+1}^*) = \Ext_R^{n+1}(B,A)$ for $n \geq 0$.
\end{proof}

While not formally dual, a similar argument using the covariant long exact sequence lets us compute \( \Ext^n_R(B,A) \) via an injective resolution of \( A \).

\begin{rmk}
Examining the proof of \cref{pro:ext-projective-resolution} shows that
in order to compute $\Ext^*_R(B, A)$,
one can replace the assumption that each $P_n$ is projective by the
assumption that $\Ext^k_R(P_n, A) = 0$ for all $k \geq 1$ and all $n \geq 0$.
(These are often called the $\Mod{R}(-, A)$-acyclic objects.)
Since we are only requiring this for one $A$, it does not suffice
to ask that $\Ext^1_R(P_n, A) = 0$, as it did in \cref{pro:HoTT-proj}.
The dual statement holds as well.
Since projective and injective resolutions may not exist in general,
this may be useful for some computations.
\end{rmk}

\subsection{Ext of finitely presented modules over (constructive) PIDs} \label{ssec:fp}
In \cref{ssec:projectives}, we will see examples which demonstrate that higher Ext groups of abelian groups do not necessarily vanish.
The main result of this section is that finitely presented abelian groups \(B\) merely have projective dimension at most \(1\), and consequently \(\Ext_{\Zb}^n(B,-)\) vanishes for \(n > 1\).
This is true more generally for finitely presented modules over \emph{principal ideal domains}, in the constructive sense of~\cite{LQ}.
When we refer to results of~\cite{LQ} in this section, we mean that these specific statements and their proofs are straightforward to translate into HoTT.
We do not make any general claims about which parts of~\cite{LQ} can be translated into HoTT.

Before turning to the constructive definition of a PID, we briefly discuss finitely presented modules.

\begin{dfn} \label{dfn:fg-fp}
  Let \(R\) be a ring and let \(A\) be an \(R\)-module.
  \begin{enumerate}[(i)]
  \item $A$ is \textbf{finitely generated} if there merely exists an epimorphism \(R^n \thra A\), for some \(n : \Nb\).
  \item $A$ is \textbf{finitely presented} if there merely exists an epimorphism \(p : R^n \thra A\), for some \(n : \Nb\), such that the kernel of \(p\) is finitely generated.
  \end{enumerate}
\end{dfn}

If \(A\) is finitely presented, then \cite[Lemma~IV.1.0]{LQ} implies that any homomorphism \(R^n \thra A\) has finitely generated kernel.
Moreover, Proposition~4.2(i) of loc.\ cit.\ says that a quotient \(A / I\), where \(I\) is a finitely generated submodule of \(A\), is also finitely presented.
These facts play a role later in this section.

We now recall the constructive definition of a PID, along with other relevant notions from~\cite{LQ}.

\begin{dfn} \label{dfn:PID}
  Let \(R\) be a commutative ring, and write \( x \mid y \defeq (\sum_{a : R} ax = y) \) for \(x, y : R\).
  \begin{enumerate}[(i)]
  \item $R$ is is an \textbf{integral domain} if every element \(x : R\) is either equal to \(0\) or \textbf{regular}: the (left) multiplication map \( y \mapsto xy : R \to R\) is a monomorphism (of \(R\)-modules).
  \item A \define{greatest common divisor} of $x, y : R$ is an
element $g$ such that the following holds:
    \[ g \mid x \times g \mid y \times \bigg(\prod_{z : R} (z \mid x \times z \mid y) \to z \mid g \bigg) .  \]
  \item $R$ is a \textbf{B\'ezout ring} if for every \(x, y : R\) there merely exist \(u, v : R\) such that \(ux + vy\) is a greatest common divisor of \(x\) and \(y\).
    The data of such a \(u\) and \(v\) is called a \textbf{B\'ezout relation} for \(x\) and \(y\).
  \item $R$ is a \textbf{B\'ezout domain} if it is both a B\'ezout ring and an integral domain.
  \item $R$ is a \textbf{principal ideal domain (PID)} if it is a B\'ezout domain, and every ascending chain of finitely generated ideals merely admits two equal consecutive terms.
  \end{enumerate}
\end{dfn}

This definition of PIDs might seem foreign to classically trained mathematicians, so we take a moment to give some context.

\begin{dfn}
  An ideal \(I\) of a ring \(R\) is \textbf{principal} if the proposition \( \Tr{ \sum_{a : R} Ra = I} \) holds.
\end{dfn}

It is not true in our setting that all ideals of \(\Zb\) are principal.
This is for a good reason: in models, there may be ``local'' ideals which have no ``global'' generators.
However, all \emph{finitely generated} ideals of \(\Zb\) are principal in our setting, and it is straightforward to verify that \(\Zb\) is a PID in the sense of \cref{dfn:PID}.
Indeed, in \(\Zb\) one can actually \emph{compute} greatest common divisors and B\'ezout relations---they don't just merely exist.
The ascending chain condition actually computes as well: using the following lemma it reduces to checking equality of principal ideals, which one can do since \(\Zb\) has decidable equality.

\begin{lem} \label{lem:bezout-fg-ideals}
  Suppose \(R\) is a B\'ezout ring.
  Any finitely generated ideal of \(R\) is principal.
\end{lem}

\begin{proof}
  The existence of B\'ezout relations means that every ideal of \(R \) that is generated by two elements is principal.
  The claim follows by recursion.
\end{proof}

The reason for the additional ``Noetherianity'' condition in our definition of PID is that it is needed to compute Smith normal forms---see~\cite[p.~209]{LQ} for further discussion.
We also get that any finitely presented module over a PID merely splits into a free part and a product of cyclic modules.
Using additivity of \(\Ext^1_{\Zb}(-,A) \) (\cref{pro:ext-additive}), projectivity of \(\Zb\), and that \(\Ext^1_{\Zb}(\Zb/n, A) \simeq A/n\)~\cite[Corollary~21]{Fla23} we deduce:

\begin{pro} \label{pro:ext-fp}
  Let \(B\) be a finitely presented abelian group, and write \( B \simeq (\bigoplus_{i=1}^k \Zb/d_i) \oplus \Zb^r\) for the decomposition which merely exists by~\cite[Prop.~IV.7.3]{LQ}.
  For any abelian group \(A\), we merely have an isomorphism
  \( \Ext^1_\Zb(B,A) \simeq \prod_{i = 1}^k A / d_i. \)  \qed
\end{pro}

From the existence of Smith normal forms over PIDs, we also deduce the following.

\begin{pro} \label{pro:PID-diagonalize}
  Suppose \(R\) is a PID.
  For any \(R\)-linear morphism \( \alpha : R^m \to R^n\), there merely exist \(R\)-linear automorphisms \(\phi\) and \(\psi\) of respectively \(R^m\) and \(R^n\) such that \(\psi \alpha \phi\) sends the \(i^{th}\) basis vector in \(R^m\) to a multiple of the \(i^{th}\) basis vector in \(R^n\) for \(1 \leq i \leq \on{min}(m,n)\).
\end{pro}

\begin{proof}
Follows from \cite[Proposition~IV.7.3(i)]{LQ}, whose proof is straightforward to carry out in HoTT.
\end{proof}

Using the proposition, we can prove the following generalization of \cref{lem:bezout-fg-ideals}.

\begin{pro} \label{pro:PID-merely-free}
  Let \(R\) be a PID, and \(n : \Nb\).
  A finitely generated submodule of \(R^n\) is merely free.
\end{pro}

\begin{proof}
  Let \(K\) be a finitely generated submodule of \(R^n\) for some \(n: \Nb\).
  We need to show that there merely exists some \(k : \Nb \) and an isomorphism \(R^k \simeq K\).
  By our assumption that \(K\) is finitely generated, there merely exists an epimorphism \(R^l \thra K\) for some \(l : \Nb\).
  Write \(\alpha : R^l \to K \to R^n\) for the composite homomorphism.
  Since we are proving a proposition, we may assume that the matrix of \(\alpha\) is diagonal, in the sense of \cref{pro:PID-diagonalize}.
  The elements on the diagonal are either regular or zero, by integrality, and we may consider the standard basis elements \(e_i : R^l\) such that \(\alpha(e_i)_i\) is regular.
  Thus we get an inclusion \(R^k \subseteq R^l\) induced by including these basis elements \(e_i\).
  Finally, the composite homomorphism \( p : R^k \to R^l \to K \) is necessarily an epimorphism, since we only threw away basis elements of \(R^l\) which are sent to \(0\) by \(\alpha\).
  By construction, the restriction of \(\alpha\) to \(R^k\) is an embedding, thus \(p\) factors an embedding and must be one itself.
  It follows that \(p\) is an isomorphism.
\end{proof}

Recall that \(\Ext^1_R(B,A)\) is itself an \(R\)-module whenever \(R\) is commutative.
For a PID \(R\) we deduce from the above that \(\Ext^1_R(B,A)\) is finitely presented (as an \(R\)-module) whenever \(A\) and \(B\) are, and moreover that \(\Ext_R^n(B,-)\) vanishes for \(n > 1\).

\begin{cor} \label{cor:ext-vanishes-fp}
  Let \(R\) be a PID.
  If \(B\) is a finitely presented \(R\)-module and \(A\) is any \(R\)-module,
  then \( \Ext_R^n(B,A) = 0 \) for \(n>1\).
  If \(A\) is also finitely presented, then so is the \(R\)-module \(\Ext_R^1(B,A)\).
\end{cor}

\begin{proof}
  Let \(B\) be a finitely presented \(R\)-module and let \(A\) be any \(R\)-module.
  Since we are proving a proposition, we may assume we have a short exact sequence
  \( K \to R^n \to B \) where the kernel \(K\) is finitely generated.
  The previous proposition lets us moreover assume that \(K\) is actually free of finite rank \(m\).
  Thus the short exact sequence gives a projective resolution of \(B\), and the claim for \(n > 1\) immediately follows by computing \(\Ext^n_R(B,A)\) using this projective resolution (\cref{pro:ext-projective-resolution}).
  The calculation of \(\Ext_R^1(B,A)\) using this projective resolution yields an exact sequence
  \[ A^m \to A^n \to \Ext_R^1(B,A) \to 0 . \]
  Hence, if \(A\) is finitely presented, then \(\Ext^1_R(B,A)\) is a quotient of the finitely presented module \(A^n\) by a finitely generated submodule (the image of \(A^m \to A^n\)), which is finitely presented.
\end{proof}

\subsection{Ext of \texorpdfstring{\(\Zb G\)}{ZG}-modules} \label{ssec:Ext-ZG}
Any construction in homotopy type theory can be carried out ``in context,'' meaning that the terms going into a particular construction may themselves depend on some extraneous variable.
In an \(\infty\)-topos, the corresponding thing is to work in a \emph{slice} of your base \(\infty\)-topos.
In this section, we work in the context of a pointed, connected type \(X\) whose base point will be denoted by \(\pt : X\).
We will see that abelian groups in the context of \(X\) correspond to modules over the group ring \(\Zb \pi_1(X)\), and we will discuss our Ext groups in this setting.

By ``an abelian group in the context of \(X\)'' we mean a family \(X \to \Ab\).
These form a category whose morphisms are families of the form \(\prod_{x : X} \Ab(A_x, B_x)\) for \(A, B : X \to \Ab\), i.e., natural transformations (the naturality squares are automatic since \(X\) is just a type).
Since \(\Ab\) is a 1-type, it is equivalent to consider families on the 1-truncation of \(X\).
The latter is equivalent to \(B \pi_1(X)\) since \(X\) is pointed and connected.
To emphasize that no truncation assumptions are needed, we prefer to work with \(X\) in this section.

We begin by constructing the group ring \(\Zb G\) for a group \(G\).
For this we use the coproduct of abelian groups, which has various constructions---see, e.g.,~\cite{Fla22, LLM}.

\begin{cst} \label{cst:zg}
  Let \(G\) be a group.
  We construct the \textbf{group ring \boldmath{\(\Zb G \)}} as follows.
  The underlying abelian group of \(\Zb G\) is the coproduct \( \bigoplus_{G} \Zb \),
  which is the free abelian group on the set \(G\).
  To define a bilinear map
  \( \Zb G \otimes_{\Zb} \Zb G \to \Zb G \)
  it suffices, by the tensor-hom adjunction and the universal property of coproducts, to give a function
  \[ G \to G \to \Zb G.\]
  For this we supply the map \((g, h) \mapsto 1_{gh}\), where \(1_{gh}\) is the unit in the \(gh\)-summand of \(\bigoplus_G \Zb\).
  The resulting binary operation on \(\Zb G\) is bilinear, by construction.
  We need to check that \(1_{e} : \Zb G\) is a two-sided unit (where \(e : G\) is the unit), and that multiplication is associative.

  Under the equivalences
  \[ (\Zb G \to \Zb G) \ \simeq \ \big(\prod_{g : G} (\Zb \to \Zb G)\big) \ \simeq ( G \to \Zb G) \]
  the identity map on \(\Zb G\) corresponds to \(g \mapsto 1_g : G \to \Zb G\).
  Since \(eg = g\) for \(g : G\), we see that \(1_e \cdot (-) : \Zb G \to \Zb G\) is the identity.
  Similarly for \((-) \cdot 1_e\).

  For associativity, simply observe that the two maps
  \[ (-) \cdot \big( - \cdot - \big),\; \big( - \cdot - \big) \cdot (-) : (\Zb G)^3 \to \Zb G \]
  both correspond to the map \((g,h,k) \mapsto 1_{ghk} : G \to G \to G \to \Zb G\), since \(G\) is associative.
\end{cst}

Before our next statement, we specify that by an \define{invertible} element of a (possible non-commutative) ring \(R\), we mean an element with a specified two-sided inverse.%
\footnote{It is equivalent to require separate left and right inverses, since one can prove that these must agree (when both exist).}
The type of two-sided inverses of a fixed element is a proposition.
We write \(\units{R}\) for the group of invertible elements of a ring \(R\).

\begin{pro}
  The group ring functor \(\Zb(-) : \Group \to \Ring \) is left adjoint to \(\units{(-)}\).
\end{pro}

\begin{proof}
  We construct a bijection
  \( \Ring(\Zb G, R) \simeq \Group(G, \units{R}) \) which is natural in \(G\) and \(R\).
  By the universal property of the coproduct, we already have a bijection
  \( \Mod{\Zb}(\Zb G, R) \simeq (G \to R) \).
  If a homomorphism on the left-hand side is in fact a ring homomorphism, then the corresponding map \(G \to R\) lands in \(\units{R}\), since ring homomorphisms are required to preserve the unit.
  Thus what we need to show is that a map \(\phi : G \to \units{R}\) is a group homomorphism if and only if the induced homomorphism \(\hat{\phi} : \Zb G \to R\) of abelian groups is a ring homomorphism.

  A map \(\phi : G \to \units{R}\) is a group homomorphism if and only if the two maps
  \[ \phi(-) \cdot \phi(-), \ \phi\big((-) \cdot (-)\big) : G^2 \to \units{R} \]
  coincide.
  This happens if and only if the two maps
  \( \hat{\phi}(-) \cdot \hat{\phi}(-), \ \hat{\phi} \big( (-) \cdot (-) \big) : \Zb G^2 \to R \)
  coincide.
  In other words, \(\phi\) is a group homomorphism if and only if \(\hat{\phi}\) is a ring homomorphism.
\end{proof}

Using the previous proposition, we relate the category \( X \to \Ab\) of abelian groups in the context of \(X\) to \(\Zb \pi_1(X)\)-modules.
Recall that for any abelian group \(M\) and ring \(R\), an \(R\)-module structure on \(M\) corresponds to a ring homomorphism \(R \to \Ab(M,M) \).

\begin{pro} \label{pro:zg-mod}
  We have an equivalence of 1-categories
  \( \Mod{\Zb \pi_1(X)} \simeq (X \to \Ab). \)
\end{pro}

\begin{proof}
  Using the previous proposition and uniqueness of deloopings of homomorphisms between groups, we have the following equivalences of types:
  \begin{align*}
    \Mod{\Zb \pi_1(X)} \ &\simeq \ \sum_{M : \Ab} \Ring\big(\Zb \pi_1(X), \Ab(M,M)\big) \\
                &\simeq \ \sum_{M : \Ab} \Group\big(\pi_1(X), \on{Aut}_{\Zb}(M)\big) \\
                &\simeq \ \sum_{M : \Ab} \big(B\pi_1(X) \pto (\Ab, M)\big) \\
                &\simeq \ (B\pi_1(X) \to \Ab) \ \simeq (X \to \Ab),
  \end{align*}
  where in the second-last line \( (\Ab, M) \) is the type \(\Ab\) equipped with the base point \(M\), and the last line uses that \(B \pi_1(X)\) is the \(1\)-truncation of \(X\).
  It is straightforward to make this association into a functor which is an equivalence of categories.
\end{proof}

Given a family \(A : X \to \Ab\), the abelian group underlying the corresponding \(\Zb \pi_1(X)\)-module is \(A_\pt\), the evaluation of \(A\) at the base point \(\pt : X\).
Accordingly, we may view \(A_\pt\) either as an abelian group or as a \(\Zb \pi_1(X)\)-module, depending on context.

An example of particular interest to us is the following.

\begin{pro} \label{pro:zg-ext-description}
  For \(B, A : X \to \Ab\), the abelian group \(\Ext_{\Zb}^n(B_\pt,A_\pt)\) is naturally a \(\Zb \pi_1(X)\)-module.
\end{pro}

\begin{proof}
  Apply \cref{pro:zg-mod} to the family \(x \mapsto \Ext_{\Zb}^n(B_x,A_x)\).
\end{proof}

When $n = 1$, we can understand the action via the following lift to $\SES_{\Zb}$.
For any \(x : X\), consider the type of short exact sequences from \(A_x\) to \(B_x\):
\[ x \longmapsto \SES_{\Zb} (B_x, A_x) \ : \ X \longrightarrow \Type. \]
This family defines a \(\Omega X\)-action on \(\SES_{\Zb}(B_*,A_*)\), and one can check that the action of an element \(g : \Omega X\) on a short exact sequence \(E\) is given by
\begin{equation} \label{eq:ses-action}
  g \cdot E \defeq (A_\pt \xra{i_E \circ g^{-1}} E_\pt \xra{g \circ p_E} B_\pt) ,
\end{equation}
where we have used the action of \(g\) on \(A_\pt\) and \(B_\pt\).
\cref{lem:pullback-pushout-iso} gives an alternative description in terms of pullbacks and pushouts.
On components, this gives the action of $\pi_1(X)$ on $\Ext_{\Zb}^1(B_*,A_*)$ from \cref{pro:zg-ext-description}.

For \(n > 1\), one gets a \(\pi_1(X)\)-action on \(\Ext^n_{\Zb}(B_\pt,A_\pt)\) which is similar to \cref{eq:ses-action}, but with a (representative of a) longer extension in place of a short exact sequence.

The following theorem identifies the type of fixed points of the action~\eqref{eq:ses-action}.

\begin{thm}\label{thm:Z-vs-ZG}
  For any \(B, A : X \to \Ab\), we have an equivalence
  \[ \prod_{x : X} \SES_{\Zb}(B_x,A_x) \ \simeq \ \SES_{\Zb \pi_1(X)}(B_\pt,A_\pt). \]
\end{thm}

\begin{proof}
  An element of \(\prod_{x:X} \SES_{\Zb}(B_x, A_x)\) is easily seen to consist of a family \(E : X \to \Ab\) along with two sections \(i : \prod_{x:X} \Mono{\Zb}(A_x, E_x)\) and \(p : \prod_{x : X} \Epi{\Zb}(E_x, B_x)\) such that the proposition
  \( \prod_{x : X} \mtt{IsExact}(i_x,p_x) \) holds.
  Under \cref{pro:zg-mod}, \(i\) corresponds to a monomorphism \(A_\pt \to E_\pt\) of \(\Zb \pi_1(X)\)-modules and \(p\) corresponds to an epimorphism \(E_\pt \to B_\pt\).
  The proposition \( \prod_{x : X} \mtt{IsExact}(i_x,p_x) \) holds if and only if it holds at the base point of \(X\), since \(X\) is connected.
  In other words, it holds if and only if \(A_\pt \to E_\pt \to B_\pt\) is an exact sequence of abelian groups (and hence of \(\Zb \pi_1(X)\)-modules).
\end{proof}

\begin{cor} \label{cor:Ext-local-cohomology}
  For any \(M\) in the abelian category \(X \to \Ab\), we have a group isomorphism
  \[ H^1(X; M) \ \simeq \ \Ext_{\Zb \pi_1(X)}^1(\Zb, M_\pt), \]
  where the left-hand side is the cohomology of \(X\) with local coefficients in \(M\), and \(\Zb\) on the right has trivial \(\Zb \pi_1(X)\)-action.
\end{cor}

\begin{proof}
  Since \(\Zb\) is projective as an abelian group, by \cref{pro:HoTT-proj} we have that \(\SES_{\Zb}(\Zb, M_x)\) is connected, for any $x : X$.
  By \cref{cor:loops-ses1}, its loop space is \(\Ab(\Zb, M_x) \simeq M_x\).
  It follows that we have an equivalence \(\SES_{\Zb}(\Zb, M_x) \simeq \K{M_x}{1}\) which is natural in \(M_x\).
  Thus we get an equivalence
  \[ \prod_{x : X} \K{M_x}{1} \ \simeq \ \prod_{x : X} \SES_{\Zb}(\Zb, M_x) \]
  which is natural in $M$.
  The set-truncation of the left-hand side is by definition \(H^1(X; M)\) (see~\cite{shulman_sss}), and the set-truncation of the right-hand side is \(\Ext_{\Zb \pi_1(X)}^1(\Zb,M_\pt)\) by the previous theorem.
  After truncating the equivalence above,  we get a natural bijection which is an isomorphism by \cref{lem:naturality-implies-homomorphism}.
\end{proof}

\section{Ext in an \texorpdfstring{\(\infty\)}{∞}-topos}
\label{sec:interpretation}

Statements in HoTT can be interpreted into an \(\infty\)-topos~\cite{dBB,dB,KL,LS,Shu19}.
In this section, we study the interpretation of the constructions and results from \cref{sec:ext}.
Our precise setup is explained in the section on foundations just below.

Results about rings and modules in HoTT apply to ring or module objects in \(\Xc\), which we stress are 0-truncated.
Accordingly, these objects live in the sub-\(\infty\)-category of 0-truncated objects in \(\Xc\), which is a 1-topos~\cite[Theorem~6.4.1.5]{HTT}.
In particular, if \(R\) is a ring object in \(\Xc\), then the category of \(R\)-modules is equivalent to a category of ordinary sheaves of modules.
Such categories have been extensively studied, and the reader may for example refer to~\cite[Chapter~18]{KS06} for background.

In \cref{ssec:SES}, we work out the interpretation \(\iSES{R}(B,A)\) of the type \(\SES_R(B,A)\) of short exact sequences, given a ring object \(R\) and two \(R\)-module objects \(A\) and \(B\) in \(\Xc\).
(Our font usage is explained below.)
The object \(\iSES{R}(B,A)\) is shown to classify short exact sequences \(A \to E \to B\) of \(R\)-modules in \(\Xc\) (\cref{pro:ses-classifies}).
From this we deduce that the set of path components of the space of global points of this object recovers the usual \emph{external} Yoneda Ext groups (\cref{cor:recover-Ext}).

Our next objective is to understand the interpretation \(\iExt{R}^n(B,A)\) of our Ext groups, which are abelian group objects in \(\Xc\).
In special cases (\cref{pro:global-points-HoTT-ext-projective-resolution,cor:global-points-sheaf-Ext}),
we show that the global points of $\iExt{R}^n(B,A)$ recover the ordinary Ext groups.
But this fails in general.
Indeed, we give examples showing that either one can vanish without the other one vanishing
(\cref{exa:external-not-internal-projective,exa:internal-not-external-projective,exa:ext-G-modules}).
However, we show that in many cases $\iExt{R}^n(B,A)$ recovers a known construction.
In any $1$-topos, one can define \emph{sheaf} Ext groups (\cref{dfn:sheaf-Ext}) by taking the right derived functors of the internal hom of modules, using the existence of enough (external) injectives.
(The name ``sheaf'' Ext is used because one often works in a category of sheaves; the name ``local'' Ext is also used.)
We can extend this to an $\infty$-topos $\Xc$, by considering sheaf Ext in the \(1\)-topos of 0-truncated objects in \(\Xc\).
When sets cover in $\Xc$ (see \cref{dfn:n-types-cover}), we show that $\iExt{R}^1(B,A)$ agrees with sheaf Ext (\cref{thm:recover-sheaf-ext}).
We do this by showing that for such $\Xc$, injectivity of modules in HoTT corresponds to internal injectivity (\cref{cor:hott-internal-injectivity}).
Since external injectives are always internally injective,
it follows that our Ext groups can also be computed using externally injective resolutions, and therefore that they agree with sheaf Ext.
A consequence of this is that in this setting our Ext groups only depend on the \(1\)-topos of 0-truncated objects in \( \Xc \) (\cref{cor:ext-determined-by-sets}).

We also study various notions of projectivity in \cref{ssec:projectives}, and provide a computation of our Ext groups using a resolution which is projective in the sense of HoTT in \cref{pro:ext-computation}.
This computation demonstrates, in particular, that our higher Ext groups need not vanish over the ring object \(\Zb\).
Lastly, \cref{ssec:Ext-BG} contains a detailed study of our Ext groups over a pointed, connected type \(X\) and over a group ring \(\Zb G\).
The considerations in this final section are meant to illustrate and exemplify the theory developed throughout \cref{sec:interpretation}, in addition to being of interest in their own right.

\addtocontents{toc}{\SkipTocEntry}
\subsection*{Foundations}
We explain our setup for interpreting HoTT into the \(\infty\)-topos \(\Xc\).
We assume an inaccessible cardinal \(\kappa\) for the entirety of \cref{sec:interpretation}.
Formally, the interpretation of HoTT lands in a type-theoretic model topos \(\cat{M}\) presenting \(\Xc\), which always exists~\cite{Shu19}.
The model topos \(\cat{M}\) admits a univalent universe which classifies relatively \(\kappa\)-presentable fibrations.
This universe allows us to interpret HoTT with a single universe.
Constructions in \(\cat{M}\) present constructions in \(\Xc\), and we are interested in studying the fruits of our labour in the latter.
We emphasize that our constructions can be built from truncations and pushouts, which are modelled in \(\infty\)-toposes; no other HITs are needed.
These constructions are all uniquely determined up to equivalence by their universal properties coming from the interpretation of the various type constructors, and this obviates the need to explicitly work with \(\cat{M}\).

Moreover, it is shown in~\cite{Stenzel} that the univalent universe in \(\cat{M}\) presents an object classifier \(\iuniv : \iTypedom \to \iType\) for relatively \(\kappa\)-compact morphisms in \(\Xc\)~\cite[Section~6.1.6]{HTT}.
This means that the mapping space \(\Xc(X,\iType)\) is naturally equivalent to the space \(\core{\Xc \kslice X}\) of relatively \(\kappa\)-compact maps into \(X\) in \(\Xc\), and lets us precisely determine the objects and structures in \(\Xc\) which are classified by the universes (of types, and of modules) that we consider.
We write \(\Xck\) for the sub-\(\infty\)-category of \(\kappa\)-compact objects in \(\Xc\).

Our results from the previous section concern truncated objects such as modules, and types of short exact sequences.
The truncation level makes the interpretation particularly straightforward, and there is not much higher coherence to manage.
For this reason---and for reasons of space and interest---we allow ourselves to state and work with the result of our interpretation directly in \(\Xc\) and not make any further mention of \(\cat{M}\).

\addtocontents{toc}{\SkipTocEntry}
\subsection*{Notation and conventions}
  We write $\Xc$ for a fixed $\infty$-topos throughout this section.
  By ``topos'' we mean \emph{Grothendieck} topos unless otherwise specified.
  Fonts are used to distinguish types in HoTT, the objects obtained by interpretation in \( \Xc \), and the classical counterparts.
  For example, \( \Ext_R^n(B,A) \) will continue to mean the Ext group in HoTT constructed in \cref{sec:ext}.
  Its interpretation \( \iExt{R}^n(B,A) \) into \( \Xc \) is written in typewriter font.
  The classical external Ext groups are written in normal font \( \eExt{R}^n(B,A) \), whereas the classical sheaf Ext groups are denoted with an underline \( \sExt{R}^n (B,A) \).
  In general, we use underlines to denote traditional constructions internal to \(\Xc\), such as the internal hom \(\ihMod{R}(A,B)\) between two \(R\)-module objects \(A\) and \(B\).
  The (external) set of \(R\)-module homomorphisms is \(\eMod{R}(A,B)\), and we use that this is equivalent to the global points \( \Gamma \ihMod{R}(A,B) \) of the internal hom.

  The 1-topos of 0-truncated objects in \(\Xc\) is denoted \(\sets{\Xc}\), and we write \(\eSet{\Xc}\) for \(\sets{\Xck}\).
  Note that the inclusion \(\sets{\Xc} \to \Xc\) respects internal homs, products, effective epimorphisms, as well as other structures and properties.
  We allow ourselves to fluently pass such structures and properties between \(\sets{\Xc}\) and \(\Xc\).
  We write \(\eAb{\Ec}\) for the (abelian) category of abelian group objects in a (possibly elementary) 1-topos \(\Ec\), and define \(\eAb{\Xc} \defeq \eAb{\sets{\Xck}}\).
  The \(\infty\)-topos of spaces is denoted \(\spaces\), and we simply write \(\eAb{}\) for \(\eAb{\spaces}\), the category of ordinary (\(\kappa\)-compact) abelian groups.
  Our abelian group and module objects in \(\Xc\) are always be assumed to be \(\kappa\)-compact.

  Base points are denoted by \(\pt\), unless another name is given.

\subsection{The object of short exact sequences}\label{ssec:SES}
Let \( R \) be a ring object in \(\Xck\), i.e., a ring object in the 1-topos \(\eSet{\Xc}\),
and write $\eMod{R}$ for the category of (\(\kappa\)-compact) $R$-modules.
Statements from HoTT about rings can be interpreted into \(\Xc\) to give statements about \(R\), and statements about $R$-modules in HoTT interpret to statements about $R$-module objects in $\Xc$.
In particular,~\cite[Theorem~4.3.4]{Fla22} shows that the category of modules (in \(\Type\)) over \(R\) interprets to an internal category \(\iMod{R}\) which represents the presheaf of \(1\)-categories
\[ X \longmapsto \eMod{(X{\times}R)} \ : \ \Xc^\op \longrightarrow \Cat, \]
where \(X \times R\) is the ring object in \(\Xc / X\) obtained by pulling back.
(Here, by ``internal category'' we mean the Rezk $(1,1)$-objects of~\cite[Def.~4.1.1]{Fla22} which represent presheaves of $1$-categories on $\Xc$, as explained in loc.\ cit.)
Thus a family of modules \(X \to \iMod{R}\) in \(\Xc\) corresponds precisely to a (relatively \(\kappa\)-compact) \( (X {\times} R) \)-module in the slice \(\Xc / X\).

For any two \(R\)-modules \( A \) and \( B \) in \( \Xck \), we interpret the type \( \SES_R(B,A) \) into \( \Xc \) to get an object \( \iSES{R}(B,A) \) of short exact sequences.
We start by describing this and the interpretation of our Ext groups, for spaces:

\begin{pro} \label{pro:recover-Ext-spaces}
  Let \(R\) be a ring object in \(\spaces\) (i.e., an ordinary ring), and let \(B\) and \(A\) be \(R\)-modules.
  \begin{enumerate}[(i)]
  \item The interpretation of \(\SES_R(B,A)\) into \(\spaces\) is equivalent to the ordinary (1-truncated) space of short exact sequences from \(A\) to \(B\);
  \item The interpretation of \(\iExt{R}^n(B,A)\) into \(\spaces\) is isomorphic to the ordinary Ext group \(\eExt{R}^n(B,A)\).
  \end{enumerate}
\end{pro}

In spaces, we will also use a slight generalization of this statement where the category of \(R\)-modules is replaced by an arbitrary abelian (univalent) category.

\begin{proof}
  Using that \(\iMod{R}\) classifies \(R\)-modules~\cite{Fla22}, it is straightforward to check that the interpretation of \(\iSES{R}(B,A)\) is equivalent to the usual space \(\eSES{R}(B,A)\) of short exact sequences.
  It follows that the interpretation of our \(\Ext_R^1(B,A)\) recovers the ordinary Yoneda Ext group \(\eExt{R}^1(B,A)\).

  For \(n > 1\), the interpretation \(\iExt{R}^n(B,A)\) into spaces recovers Yoneda's definition of \(\eExt{R}^n(B,A)\) as a quotient of the space of length-\(n\) exact sequences, which is well-known to give the usual Ext groups defined using resolutions.
\end{proof}

Our present goal is to relate the object \( \iSES{R}(B,A) \) in \(\Xc\) to the external space \( \eSES{R}(B,A) \), for ring and module objects in \(\Xc\).
To do this, we require a lemma which characterizes the interpretation of the objects of epimorphisms and mono\-mor\-phisms from HoTT.

\begin{lem} \label{lem:object-of-epis}
  Let \(A\) and \(B\) be \(R\)-modules in \(\Xc\).
  The object \(\iEpi{R}(B,A)\) resulting from interpretation classifies \(R\)-module epimorphisms \(B \to A\) in \(\Xc\).
  Likewise, the object \(\iMono{R}(B,A)\) classifies \(R\)-module monomorphisms.
\end{lem}

\begin{proof}
Let us first recall that the internal hom \( \ihMod{R}(B,A) \) classifies \(R\)-module homomorphisms, meaning that we have natural equivalences
\[ \Xc(X, \ihMod{R}(B,A)) \ \simeq \ \eMod{(X{\times}R)}(X \times B, X \times A) \]
for \(X \in \Xc\).
We will restrict these equivalences to epimorphisms and monomorphisms to obtain our statement.

  Under this equivalence, a map \(f : X \to \iEpi{R}(B,A)\) corresponds to an \((X{\times}R)\)-module homomorphism \(f' : X \times B \to X \times A\) in \(\Xc / X\) that satisfies the interpretation of being \((-1)\)-connected from HoTT.
  Since \((-1)\)-connected maps in HoTT interpret to \((-1)\)-connected maps in an \(\infty\)-topos, we have that \(f'\) is a \((-1)\)-connected map over \(X\).
  This means that \(f'\) is an (effective) epimorphism between 0-truncated objects, hence also between module objects, as desired.

  The statement for \(\iMono{R}(B,A)\) is shown similarly, but using that \((-1)\)-truncated maps in HoTT correspond to \((-1)\)-truncated maps in \(\Xc\), which are monomorphisms between 0-truncated objects.
\end{proof}

We use this lemma in the proof of the following proposition, which says that the object of short exact sequences from HoTT classifies short exact sequences in \(\Xc\).
Recall that base change functors are exact and therefore preserve ring and module objects, as well as exact sequences of the latter.
This means that any morphism \(f : X \to Y\) in \(\Xc\) induces a map
\[ f^* \ : \ \eSES{(Y{\times}R)}(Y \times B, Y \times A) \longrightarrow \eSES{(X{\times}R)}(X \times B, X \times A) \]
by base change, for any two \(R\)-modules \(A\) and \(B\) in \(\Xc\).

\begin{pro} \label{pro:ses-classifies}
  Let \( A \) and \(B \) be \(R\)-modules in \( \Xck \).
  The object \( \iSES{R}(B,A) \) represents the presheaf
  \[ X \longmapsto \eSES{(X{\times}R)}(X \times B, X \times A) \ : \ \Xc^\op \longrightarrow \spaces. \]
  In particular, the (1-truncated) space \( \eSES{R}(B,A) \) is equivalent to the global points of the object \(\iSES{R}(B,A)\).
\end{pro}

\begin{proof}
  Let \( X \in \Xc \).
  Our goal is to produce equivalences of spaces
  \[ \Xc \big(X, \iSES{R}(B,A) \big) \ \simeq \ \eSES{(X{\times}R)}(X \times B, X \times A) \]
  which are natural in \(X \in \Xc\).
  Using the adjunction \(\Sigma_{X} \dashv X \times (-)\) and base-change stability of interpretation, we may replace the left-hand side above via the following natural equivalences:
  \[ \Xc \big( X, \iSES{R}(B,A) \big)
    \ \simeq \ \Xc / X \big( \id_X, X \times \iSES{R}(B,A) \big )
    \ \simeq \ \Xc / X \big( \id_X, \iSES{(X \times R)}(X \times B, X \times A) \big). \]
  The rightmost space is the global points of the object \(\iSES{(X \times R)}(X \times B, X \times A)\).
  It therefore suffices to consider the case \(X = 1\), and to construct, for any \(\infty\)-topos \(\Xc\), an equivalence
  \[ \Xc \big(1, \iSES{R}(B,A) \big) \ \simeq \ \eSES{R}(B,A) \]
  that is pullback-stable along any \(Y \to 1\) (which implies naturality of the required equivalence).

  The right-hand side above is the domain of a \((-1)\)-truncated map into the (1-truncated) space \(G\) consisting of \(\kappa\)-compact \(R\)-modules \(E\) equipped with a monomorphism \(i : A \to E\) and an epimorphism \(p : E \to B\).
  The object \(\iSES{R}(B,A)\) is the domain of a \((-1)\)-truncated map into the corresponding object \(G'\) of such things in \(\Xc\).
  By \cref{pro:recover-Ext-spaces}(i) applied to the abelian category \(\eMod{R}\), the previous lemma, and~\cite[Theorem~4.3.4]{Fla22}, we have a pullback-stable equivalence between the space of global points of \(G'\) and \(G\).
  It follows that both sides of the equivalence above are naturally fibred over \(G\), and we can therefore obtain the desired equivalence from a fibrewise bi-implication (which yields an equivalence since the maps are \((-1)\)-truncated).
  Note that pullback-stability of the equivalence \(\Xc(1,G') \simeq G\) will automatically imply pullback-stability of this equivalence between \((-1)\)-truncated maps.
  It therefore remains to check that the internal proposition \(\mtt{IsExact}(i,p)\) holds if and only \(i\) and \(p\) define an exact complex in the usual sense.

  The proposition \( \mtt{IsExact}(i,p) \) consists of a witness that the internally induced map \( A \to \mtt{ker}(p) \) is $(-1)$-connected.
  The module \( \mtt{ker}(p) \) is clearly equivalent to the externally defined kernel \( \on{ker}(p) \), both being given by the fibre over the global point \( 0 : 1 \to B \).
  Under this equivalence, the aforementioned witness implies that the induced map \( A \to \on{ker}(p) \) is surjective (i.e., \((-1)\)-connected), and vice-versa.

  In conclusion, \( \Xc \big(1, \iSES{R}(B,A) \big) \) and \( \eSES{R}(B,A) \) are naturally fibrewise equivalent as spaces over \(G\).
  This yields a pullback-stable equivalence on total spaces, as desired.
\end{proof}

Recall that \( \pi_0\eSES{R}(B,A) \) is the usual definition of the \emph{Yoneda Ext groups} (see, e.g.,~\cite{Mac63}), which recover the Ext groups defined in terms of resolutions.
Thus we have the following: 

\begin{cor}\label{cor:recover-Ext}
  We have a natural isomorphism \( \pi_0 \big(\Xc \big( 1, \iSES{R}(B,A) \big)\big) \simeq \eExt{R}^1(B,A) \) of ordinary abelian groups, for any \( R \)-modules \( A \) and \( B \) in \( \Xck \).\qed
\end{cor}

Since we do not have a good description of the (untruncated) type of length-$n$ exact sequences,
we do not have a corresponding statement for the higher Ext groups.

Note that taking global points and taking components do not commute, and it is
important for the above result that we take global points before taking components.
If we reverse the order, we get the claim that
\( \Xc \big( 1, \iExt{R}^1(B,A) \big) \simeq \eExt{R}^1(B,A) \).
We show in \cref{exa:internal-not-external-projective,exa:ext-G-modules} that this is
false in general.
However, we will see in \cref{pro:global-points-HoTT-ext-projective-resolution,cor:global-points-sheaf-Ext}
that there are situations in which the global points of $\iExt{R}^n$ agree with $\eExt{R}^n$ for all $n$.

\subsection{Comparing various notions of projectivity}
\label{ssec:projectives}

It is well-known that ordinary Ext groups of (say) modules can be computed using projective resolutions, whenever one is at hand.
In \cref{ssec:ext-proj-res} we showed that the same thing holds for the Ext groups we defined in HoTT.
Accordingly, we can compute our internal Ext groups in \(\Xc\) using resolutions which consist of modules that satisfy the interpretation of the predicate \(\mathsf{IsProjective}\) from \cref{dfn:HoTT-proj}.
An example of such a computation is given in \cref{pro:ext-computation}.
In addition, we compare internal projectivity to ordinary (external) projectivity.
There are no implications either way in general, which we demonstrate through \cref{exa:internal-not-external-projective,exa:external-not-internal-projective}.

\begin{dfn}
  Let \(R\) be a ring object in \(\Xc\).
  \begin{enumerate}[(i)]
  \item An \(R\)-module \(P\) is \define{(externally) projective} if for every epimorphism \(e : A \thra B\) in \(\eMod{R}\), the map \(e_* : \eMod{R}(P,A) \to \eMod{R}(P,B)\) of ordinary sets (or abelian groups) is an epimorphism;
  \item An \( R \)-module \( P \) is \define{internally projective} if for every epimorphism \( e : A \thra B \) in \( \eMod{R} \), the map \( e_* : \ihMod{R}(P,A) \to \ihMod{R}(P,B) \) in \( \eAb{\Xc} \) is an epimorphism;
  \item An \(R\)-module \(P\) is \define{HoTT-projective} if the interpretation of the proposition \(\mathsf{IsProjective}(P)\) from \cref{dfn:HoTT-proj} holds.
  \end{enumerate}
\end{dfn}

The external and internal notions are the usual ones which pertain to modules in a 1-topos, which for us is the 1-topos \(\sets{\Xc}\).
However, in an \(\infty\)-topos we also have the third notion of HoTT-projectivity resulting from interpretation.
We mention, to be concrete, that if \(\Xc\) is the sheaf \(\infty\)-topos on some 1-site, then \(\sets{\Xc}\) is the category of ordinary set-valued sheaves on the same site.
In this situation, ring and module objects are ordinary sheaves of rings and modules.

In general, when we say that the interpretation of a statement in HoTT ``holds''  we mean that the resulting object of \(\Xc\) has a global point.
If the statement is a proposition, then this means that the object is terminal.
Our first objective is to make a useful reformulation of HoTT-projectivity.

\begin{pro} \label{pro:hott-projectivity}
An \(R\)-module \(P\) is HoTT-projective if and only if the \((X{\times}R)\)-module \(X \times P\) is internally projective in \(\eMod{(X \times R)}\) for all \(X \in \Xc\).
\end{pro}

\begin{proof}
  Let \(P\) be an \(R\)-module in \(\Xc\).
  According to \cref{dfn:HoTT-proj}, we have
  \[ \mathsf{IsProjective}(P) \ \defeq \ \prod_{A : \Mod{R}} \prod_{B : \Mod{R}} \prod_{e : \Epi{R}(A,B)} \mathsf{IsEpi}\big(e_* : \Mod{R}(P,A) \to \Mod{R}(P,B)\big) .\]
  Interpreting \(\mathsf{IsProjective}(P)\), we get an object of \(\Xc\).
  It has a global point if and only if the projection
  \[ Q \ : \ \sum_{A, B : \Mod{R}} \ \sum_{e : \ttt{Epi}_R(A,B)} \ttt{IsEpi}(e_*) \ \longrightarrow \ \sum_{A, B : \iMod{R}} \ttt{Epi}_R(A,B) \]
  admits a section.
  This map admits a section if and only if for every map \(f : X \to \sum_{A, B : \Mod{R}} \ttt{Epi}_R(A,B) \) there is a section of the pullback \(f^*(Q) \in \Xc / X\), since we can take \(f\) to be the identity map.
  Such a map \(f\) is equivalent to the data of two \((X{\times}R)\)-modules \(A\) and \(B\) over \(X\) along with an epimorphism \(e : A \to B\).
  Here we have used~\cite[Theorem~4.3.4]{Fla22} which says that \(\iMod{R}\) classifies module objects, and \cref{lem:object-of-epis}.
  By definition, we have that \(f^*(Q) = \ttt{IsEpi}(e_*)\), where
  \(e_* : \ihMod{(X{\times}R)}(X \times P,A) \to \ihMod{(X{\times}R)}(X \times P, B) \)
  is the post-composition map.
  This proposition \(f^*(Q)\) holds if and only if \(e_*\) is an epimorphism.

  In summary, the statement \(\ttt{IsProjective}(P)\) holds if and only if for every \(X \in \Xc\), all \((X{\times}R)\)-modules \(A\) and \(B\), and every \(R\)-module epimorphism \(e : A \to B\), the aforementioned post-composition map \(e_*\) is an epimorphism.
  But this is exactly the statement that \(X \times P\) is an internally projective \((X{\times}R)\)-module for every \(X \in \Xc\).
\end{proof}

Clearly, HoTT-projectivity always implies internal projectivity.
The converse holds for \(\infty\)-toposes in which internal projectivity of modules is stable by base change.
We do not know whether this is always true%
\footnote{David Wärn has now shown that internal projectivity does \emph{not} imply HoTT-projectivity; see \cite[Theorem~7]{Warn24}.}%
, but it is true for spaces, as we show in \cref{pro:hott-projectivity-spaces}.

Next we show that certain free modules are HoTT-projective.
Our proof uses the following lemma, due to Alex Simpson for internal projectivity of \emph{objects}
in a \emph{1}-topos~\cite{Simpson13}, written up on the nLab\footnote{See \href{https://ncatlab.org/nlab/revision/internally+projective+object/13}{internally projective object (rev.~13)} on the nLab.}.
(See also \cite[Lemma~11]{Warn24}.)
The definition of internal projectivity of objects is the same as for modules, but using the internal hom of objects.
It is straightforward to check that Simpson's proof goes through for objects of an \(\infty\)-topos as well, providing us with:

\begin{lem}\label{lem:internal-projectivity-objects-stable}
  Let \(P \in \Xc\) be an internally projective object.
  Then \(X \times P\) is an internally projective object in \(\Xc / X\) for all \(X \in \Xc\). \qed
\end{lem}

Combining \cref{lem:internal-projectivity-objects-stable} with an argument similar
to that of \cref{pro:hott-projectivity}, we deduce:

\begin{pro}\label{pro:projective-objects}
  An object \(P \in \Xc\) is internally projective if and only if it
  satisfies the interpretation of being a projective object in HoTT. \qed
\end{pro}

Given a \(0\)-truncated object \(S\) in \(\Xc\), we can form the free \(R\)-module \(R(S)\) on this object, for any ring object \(R\).
This free \(R\)-module is the interpretation of the free \(R\)-module on a set in HoTT.

\begin{pro}\label{pro:hott-projective-free}
  Let \(R\) be a ring object in \(\Xc\), and let \(P\) be a \(0\)-truncated, internally projective object.
  The free \(R\)-module \(R(P)\) on \(P\) is HoTT-projective.
\end{pro}

\begin{proof}
  By \cref{pro:projective-objects}, \(P\) satisfies the interpretation of being a projective set in HoTT.
  The free \(R\)-module \(R(P)\) on a projective set is projective in HoTT, so we are done.
\end{proof}

We use this proposition to compute an example of our Ext groups \(\Ext_R^n\) in \cref{pro:ext-computation}.
Before turning to this example, we observe that internal projectivity (and thus HoTT-projectivity) implies external projectivity in certain situations.
This has some interesting consequences.

\begin{pro}\label{pro:internal-implies-external-projectivity}
  Let \(\Ec\) be a (possibly elementary) 1-topos, equipped with a ring object \(R\).
  If the global points functor \(\Gamma : \Ec \to \Set\) preserves epimorphisms, then internal projectivity of \(R\)-modules implies external projectivity.
\end{pro}

\begin{proof}
  The statement easily follows by identifying the external hom of \(R\)-modules as the global points of the corresponding internal hom, and then using the assumption on \(\Gamma\).
\end{proof}

The previous proposition applies, for example, to any topos of presheaves on a category with a terminal object \(1\).
In that case \(\Gamma\) is represented by evaluation at \(1\), which respects both limits and colimits of presheaves.

\begin{pro}\label{pro:global-points-HoTT-ext-projective-resolution}
  Let \(R\) be a ring object in \(\Xc\), and consider two \(R\)-modules \(B\) and \(A\).
  Suppose that \(\Gamma : \eSet{\Xc} \to \Set \) preserves epimorphisms.
  If $B$ has a HoTT-projective resolution \(P_\bullet\), then we get an isomorphism
  \( \Gamma \iExt{R}^n(B, A) \ \simeq \ \eExt{R}^n(B,A) \).
\end{pro}

One can check that our assumption on \(\Gamma\) holds if and only if the induced functor \(\Gamma : \eAb{\Xc} \to \eAb{}\) is exact, and it is this latter condition that we use in the proof.

\begin{proof}
  By the interpretation of \cref{pro:ext-projective-resolution}, we can compute \(\iExt{R}^n(B,A)\) using the HoTT-projective resolution \(P_\bullet\).
  Specifically, taking internal homs we get a complex
  \begin{equation}\label{eq:int-ext-projective-resolution}
    \cdots \to \ihMod{R}(P_{n-1},A) \to \ihMod{R}(P_n,A) \to \ihMod{R}(P_{n+1},A) \to \cdots
  \end{equation}
  of abelian groups in \(\Xc\), and we have isomorphisms
  \( H^n(P_\bullet; A) \simeq \iExt{R}^n(B,A) \)
  where the left-hand side is the cohomology of the above complex in \(\eAb{\Xc}\).

  Now, by our assumption on \(\Gamma\), the previous proposition tells us that \(P_\bullet\) is an externally projective resolution of \(B\) (since HoTT-projective always implies internally projective).
  Thus we may also compute \(\eExt{R}^n(B,A)\) using \(P_\bullet\), which amounts to taking the cohomology (in \(\eAb{}\)) of the global points of the complex (\ref{eq:int-ext-projective-resolution}) above.
  Since \(\Gamma : \eAb{\Xc} \to \eAb{}\) is exact, it commutes with taking cohomology, and we therefore obtain the desired isomorphism.
\end{proof}

In the presence of enough HoTT-projectives, we deduce:

\begin{cor}\label{cor:global-points-HoTT-ext-enough-internal-projectives}
  Let \(R\) be a ring object in \(\Xc\), and suppose that \(\eMod{R}\) has enough HoTT-projectives.
  If \(\Gamma : \eSet{\Xc} \to \Set \) preserves epimorphisms, then we have natural isomorphisms
  \[ \Gamma \iExt{R}^n(B, A) \ \simeq \ \eExt{R}^n(B,A) \]
  for any two \(R\)-modules \(B\) and \(A\).
  \qed
\end{cor}

We now turn to our computation of a non-trivial \(\iExt{\Zb}^2\) in the Sierpi\'nski \(\infty\)-topos \(\Xc\) using a HoTT-projective resolution.
Since \(\Xc\) is the \(\infty\)-topos of presheaves on the arrow category \(0 \to 1\), an abelian group in \(\Xc\) consists of a homomorphism \(A_0 \leftarrow A_1\) between two ordinary abelian groups.
We write \(\yo\) for the Yoneda embedding.

Note that \(\yo(0)\) is an internally projective object in this \(\infty\)-topos, since it represents the functor sending a presheaf \( F_0 \leftarrow F_1 \) to the presheaf \(F_0 \xla{\id} F_0\), which preserves epimorphisms.
Accordingly, the corresponding free abelian group \(\Zb \yo(0)\) is HoTT-projective, by \cref{pro:hott-projective-free}.

\begin{pro} \label{pro:ext-computation}
  Consider the abelian group \(B \defeq (0 \leftarrow \Zb / 2)\) in the Sierpi\'nski \(\infty\)-topos.
  We have
  \[ \iExt{\Zb}^2(B, \Zb y(0)) \ \simeq \ B . \]
\end{pro}

\begin{proof}
  We will compute this internal Ext group using a HoTT-projective resolution of \(B\), as justified by the interpretation of \cref{pro:ext-projective-resolution}.
  Drawing objects of \(\Xc\) vertically and morphisms horizontally, the following is such a resolution \(P_*(B)\):
  \[ \begin{tikzcd}[column sep=2.5em]
      0 \dar \rar & \Zb \rar["2"] \dar["{(0,1)}"] & \Zb \rar[two heads] \dar[equals] & \Zb / 2 \dar \\
      \Zb \rar["{(2, -1)}"] & \Zb \oplus \Zb \rar["1 + 2"] & \Zb \rar & 0 \period
  \end{tikzcd} \]
  Here, \(P_0(B) = \Zb \yo(1)\) is the integer object in \(\Xc\), which is always HoTT-projective; \(P_2(B) = \Zb \yo(0)\) is HoTT-projective by the discussion just above; and \(P_1(B)\) is \(\Zb \yo(0) \oplus \Zb \yo(1)\), a direct sum of HoTT-projectives, which is HoTT-projective.
  Thus \(\iExt{\Zb}^2(B, \Zb \yo(0)) \simeq H^2(P_*(B), \Zb \yo(0)) \).
  Since $\Zb \yo(1)$ is the integer object, we have
  $\ihMod{\Zb}(\Zb \yo(1), \Zb \yo(0)) \simeq \Zb \yo(0)$,
  and one can check that $\ihMod{\Zb}(\Zb \yo(0), \Zb \yo(0)) \simeq \Zb \yo(1)$.
  Applying $\ihMod{\Zb}(-, \Zb \yo(0))$ to the first three columns above gives the cochain complex
  \[ \begin{tikzcd}[column sep=2.5em]
      \Zb \dar[equals] & \Zb \lar["2",swap] \dar["{(1,0)}"] & 0 \lar \dar \\
      \Zb & \Zb \oplus \Zb \lar["2 - 1",swap] & \Zb \lar["{(1, 2)}",swap] \period
  \end{tikzcd} \]
  The desired group $H^2(P_*(B), \Zb \yo(0))$ is the cokernel of the leftmost map,
  which is \(B = (0 \leftarrow \Zb / 2)\).
\end{proof}

\begin{rmk}
The category \((0 \to 1)\) has a terminal object (the object \(1\)), which implies that the global points functor of the Sierpi\'nski \(\infty\)-topos preserves epimorphisms between abelian group objects.
Specifically, the global points of an object \(A_0 \leftarrow A_1\) is simply \(A_1\).
From \cref{pro:global-points-HoTT-ext-projective-resolution} and the computation in the previous proposition, we deduce that \(\eExt{\Zb}^2(B, \Zb\yo(0)) \simeq \Zb / 2\) for the external Ext group.
\end{rmk}

\begin{rmk}\label{rmk:not-a-product}
  The first part of the resolution above gives rise to a short exact sequence
  \[ \begin{tikzcd}
      \Zb \rar["2"] \dar["2"] & \Zb \rar[two heads] \dar[equals] & \Zb / 2 \dar \\
      \Zb \rar[equals] & \Zb \rar & 0 \period
  \end{tikzcd} \]
  The object in the centre is clearly not the product of the kernel and the cokernel,
  even ignoring the group structures.
  In the Sierpi\'nski \(\infty\)-topos, an object is merely inhabited if and only if
  it is inhabited, so it follows that it does not merely hold that the central object
  is the product of the kernel and the cokernel.
  So while it is true that the type of length-$1$ extensions is essentially small, this cannot
  be proved by assuming that the underlying type of the middle object is the product of
  the other two types.  (See \cref{rmk:not-a-product1} for context.)
\end{rmk}

We conclude this section by studying the relation between internal and external projectivity in general.
\cref{exa:internal-not-external-projective} gives an internally projective module which is not externally projective.
Here we give an example of an externally projective abelian group that fails to be internally projective, which is an additive version of an example due to Todd Trimble.\footnote{See \href{https://ncatlab.org/nlab/revision/presentation+axiom/46\#counter}{presentation axiom (rev. 46)} on the nLab.}

Consider the poset \( \cat{C} \defeq \mathbb{N} * \{ a , b \} \), where \( a \) and \( b \) are greater than all \( n \in \mathbb{N} \),
  \[ \begin{tikzcd}[column sep=small, row sep=tiny]
    & & & & & a \\
    0 \rar & 1 \rar & 2 \rar & \cdots \ar[urr] \ar[drr] \\
    & & & & & b
  \end{tikzcd}  \]
and let \( \cat{X} \) be the \(\infty \)-topos of presheaves on \( \cat{C} \).
As above, we write \( \yo : \cat{C} \to \cat{X} \) for the Yoneda embedding.
The functor \( \Zb(-) : \eSet{\Xc} \to \eMod{\Zb} \) constructs the free abelian group on a 0-truncated object in \(\Xc\).
In particular, we may depict \( \Zb \yo(a) \) as follows:
  \[ \begin{tikzcd}[column sep=small, row sep=tiny]
    & & & & & \Zb \ar[dll, equals] \\
    \Zb  & \lar[equals] \Zb & \lar[equals] \Zb & \lar[equals] \cdots \\
    & & & & & 0 \ar[ull, hook]
  \end{tikzcd}  \]

The integer object \(\Zb\) over \(\cat{C}\) is simply the constant presheaf on the ordinary integers.

\begin{exa} \label{exa:external-not-internal-projective}
  \( \Zb \yo(a) \) is externally projective but not internally projective in \( \eMod{\Zb} \).
  It follows from \cref{pro:HoTT-proj} that there exists an $A$ in $\eMod{\Zb}$ so that
  $\eExt{\Zb}^1(\Zb \yo(a), A) = 0$ and $\iExt{\Zb}^1(\Zb \yo(a), A) \neq 0$.
\end{exa}

\begin{proof}
It is immediate that \( \Zb \yo(a) \) is externally projective, since it represents evaluation at \( a \), which preserves epimorphisms of presheaves of modules.
To show that \( \Zb \yo(a) \) is not internally projective, we construct an epimorphism \( \sigma : F \thra G \) that isn't preserved by \( \ihMod{\Zb}(\Zb \yo(a), -) \).

Let \( F : \eMod{\Zb} \) be defined as follows, with indexwise inclusions:
  \[ \begin{tikzcd}[column sep=small, row sep=tiny]
    & & & & & 0 \ar[dll, hook] \\
    \bigoplus_{n \in \mathbb{N}} \Zb  & \lar[hook] \bigoplus_{n \geq 1} \Zb & \lar[hook] \bigoplus_{n \geq 2} \Zb & \lar[hook] \cdots \\
    & & & & & 0 \ar[ull, hook] \period
  \end{tikzcd}  \]
Define \( G(n) \defeq \Zb \) for \( n \in \mathbb{N} \) and \( G(a) := 0 =: G(b) \), with all maps between natural numbers inducing identity maps.
Then we have an epimorphism \( \sigma : F \thra G \) given by addition at \( n \in \mathbb{N} \) and identities (zero maps) at \( a \) and \(b \).
However \( \eMod{\Zb}(G,F) = 0 \) since any such homomorphism must factor through \( \on{lim}_n F (n) = 0 \) on the \( \mathbb{N} \)-part of \( \cat{C} \) (and is necessarily $0$ on \( a \) and \( b \)).

Using the tensor-hom (see, e.g., \cite[Section~17.22]{stacks}) and free-forgetful adjunctions, one can check that
\( G = \Zb ( \yo(a) \times \yo(b) ) \) and \( \Zb \yo(a) \otimes_\Zb \Zb \yo(b) \) satisfy the same universal property and are therefore naturally isomorphic.
Here \(\otimes_\Zb\) is the tensor product of presheaves of modules, which is pointwise (see, e.g., Section~6.6 of loc.\ cit.).
Using the tensor-hom adjunction again, we obtain isomorphisms
\[ \ihMod{\Zb}(\Zb \yo(a), F)(b)
  \ \simeq \ \eMod{\Zb}(\Zb \yo(a) \otimes_\Zb \Zb \yo(b), F)
  \ \simeq \ \eMod{\Zb} (G ,F )
  \ \simeq \ 0. \]
On the other hand,
\( \ihMod{\Zb}(\Zb \yo(a), G)(b) = \eMod{\Zb}(\Zb \yo(a) \otimes_\Zb \Zb \yo(b), G) = \eMod{\Zb}(G,G) \) contains at least two elements: \( 0 \) and \( \id_G \).
This means that \( \sigma_* : \ihMod{\Zb}(\Zb \yo(a), F) \to \ihMod{\Zb}(\Zb \yo(a), G) \) cannot be an epimorphism, since it isn't one at \( b \).
\end{proof}

\subsection{Internal injectivity and sheaf Ext}\label{ssec:sheaf-Ext}
The goal of this section is to show that in certain \( \infty \)-toposes, the interpretation of our Ext groups from HoTT recover the \emph{sheaf} Ext groups arising in algebraic geometry.
Since sheaf Ext groups are defined using external injective resolutions (see \cref{dfn:sheaf-Ext}),
we will need to understand how these compare to the notion of injectivity
coming from HoTT.
We will make this comparison by going through a third notion of injectivity,
namely internal injectivity.

\begin{dfn}
  Let \(R\) be a ring object in \(\Xck\).
\begin{enumerate}[(i)]
\item  An $R$-module $I$ is \define{(externally) injective} if for every $R$-module monomorphism $m : A \to B$, the homomorphism $m^* : \eMod{R}(B,I) \to \eMod{R}(A,I)$ in $\eAb{}$ is an epimorphism;
\item An $R$-module $I$ is \define{internally injective} if for every $R$-module monomorphism $m : A \to B$, the homomorphism $m^* : \ihMod{R}(B,I) \to \ihMod{R}(A,I)$ in $\eAb{\Xc}$ is an epimorphism;
\item An $R$-module $I$ is \define{HoTT-injective} if it satisfies the interpretation of the proposition \( \mathsf{IsInjective}(I) \) from \cref{dfn:HoTT-inj} in $\Xc$.
\end{enumerate}
\end{dfn}

As in the projective case, injectivity and internal injectivity are the familiar notions from the 1-topos \(\eSet{\Xc}\).
In surprising contrast to the projective case we considered above, external injectivity always implies internal injectivity in a 1-topos.
This theorem is due to~\cite{Har83} for abelian groups and~\cite[Theorem~3.8]{Ble18} for modules.
The converse holds in any localic 1-topos (as Blechschmidt shows), however not every internally injective module is externally injective in general.
For example, in~\cite[pp.~259]{Har82-effacements} it is shown that the \(\Zb / 2\)-module \( \Qb / \Zb \) with trivial action is not externally injective, though it is internally injective (as an abelian group over \(B \Zb/2\)) by~\cite[Proposition~1.2(i)]{Har81} (see also \cref{pro:internally-injective-G-module} below).

\begin{rmk}\label{rmk:enough-internal-injectives}
  Let \(R\) be a ring object in \(\Xck\).
  The category \( \eMod{R} \) is equivalent to a category of modules in a 1-topos~\cite[Theorem~6.4.1.5]{HTT}, and is therefore Grothendieck abelian~\cite[Theorem~18.1.6]{KS06}.
  Consequently, it has enough external injectives.
  Since external injective are internally injective by~\cite[Theorem~3.8]{Ble18},
there are also enough internal injectives in \(\eMod{R}\).
\end{rmk}

To relate HoTT-injectivity to internal injectivity we proceed as we did in \cref{ssec:projectives} for projectivity.
A proof similar to the one of \cref{pro:hott-projectivity} gives us the following:

\begin{pro}
  An \( R \)-module \( I \) in \( \Xc \) is HoTT-injective if and only if
  the \( (X{\times}R) \)-module \( X \times I \) is internally injective in
  \( \Xc \kslice X \) for all \( X \in \Xc \).
  \qed
\end{pro}

Clearly, every HoTT-injective module is internally injective.
If we could show that internal injectivity is stable by base change in an \(\infty\)-topos, then the two notions would coincide.
However, we do not know whether internal injectivity is stable by base change in a general \(\infty\)-topos.
We will show below that it holds in certain situations.

In~\cite{Har83}, Harting showed that internal injectivity is stable by base change for abelian groups in any elementary \(1\)-topos.
The same holds for modules (see the discussion immediately after~\cite[Proposition~3.7]{Ble18}).
It follows that the same is true in an $\infty$-topos for base change
by a $0$-truncated object:

\begin{lem} \label{lem:internal-injectivity-stable-sets}
  Let \( X \in \Xc \) be \(0\)-truncated.
  Base change \( X \times (-) : \eMod{R} \to \eMod{(X{\times}R)} \) over \( X \) preserves internal injectivity.\qed
\end{lem}

A key fact in the converse direction is that internal injectivity \emph{descends along effective epimorphisms}.
This is essentially a corollary of the fact that base change along effective epimorphisms reflects effective epimorphisms (more generally, connected maps)~\cite[Proposition~6.5.1.16(6)]{HTT}.

\begin{lem} \label{lem:internal-injectivity-descends}
  Let \(V\) be a \((-1)\)-connected object of \(\Xc\), and let \(I\) be an \(R\)-module.
  If \(V \times I\) is internally injective as a \((V{\times}R)\)-module over \(V\), then \(I\) is internally injective.
\end{lem}

The same result holds for internal projectivity as well, with only minor changes to the proof.

\begin{proof}
  Suppose \(V \times I\) is internally injective as a \((V{\times}R)\)-module, and let \( i : A \hra B\) be a monomorphism of \(R\)-modules.
  Using that base change preserves internal homs, consider the following diagram:
  \[ \begin{tikzcd}[column sep=large]
      \ihMod{(V{\times}R)}(V \times B, V \times I) \rar[two heads] \dar["{(V \times i)^*}"] \ar[dr, phantom, "\lrcorner" at start] & \ihMod{R}(B, I) \dar["i^*"] \\
      \ihMod{(V{\times}R)}(V \times A, V \times I) \dar \rar[two heads] \ar[dr, phantom, "\lrcorner" at start] & \ihMod{R}(A, I) \dar \\
      V \rar[two heads] & 1 \period
    \end{tikzcd} \]
  Here the two-headed arrows signify effective epimorphisms (which are pullback-stable).
  By assumption, \(V \times I\) is internally injective, so \((V \times i)^*\) is an effective epimorphism (since monomorphisms are stable by base change).
  Thus \(i^*\) factors an effective epimorphism, and must therefore be one itself.
  We conclude that \(I\) is internally injective, as desired.
\end{proof}

One can also prove the previous lemma by using that pullback along an
effective epimorphism is conservative~\cite[Lemma~6.2.3.16]{HTT}.
We apply this to the map $\ttt{IsEpi}(i^*) \to 1$, which is an equivalence
precisely when $i^*$ is an epimorphism.

We now introduce conditions on \(\Xc\) which will imply that internal injectivity is stable by base change.

\begin{dfn} \label{dfn:n-types-cover}
  Let \(n \geq -1\) be a truncation level.
  An object \(X \in \Xc\) is \define{covered by an \(n\)-type} if there exists an \(n\)-type \(V\) along with an effective epimorphism \(V \thra X\).
  If all objects of \(\Xc\) are covered by \(n\)-types, then \define{\(n\)-types cover} in \(\Xc\).
  When \(n = 0\), we say that \define{sets cover}.
\end{dfn}

Sets cover in any \(\infty\)-topos of \(\infty\)-sheaves on a \(1\)-category, since any such sheaf can be covered by a coproduct of representables.

Note that the condition that $n$-types cover in $\Xc$ is not the interpretation of
the corresponding concept from HoTT.
(See~\cite[Exercise~7.9]{hottbook}, \cite[Definition~5.2]{Chr21}
and the nLab\footnote{See \href{https://ncatlab.org/nlab/revision/n-types+cover/6\#in_homotopy_type_theory}{$n$-types cover (rev.~6)} on the nLab.}.)
The HoTT notion only requires that \(V\) and the effective epimorphism \emph{merely} exist.
On the other hand, it requires this in every slice.

By combining the two previous lemmas with this definition above, we obtain the following result.

\begin{pro} \label{pro:internal-injectivity-covered-by-a-set}
  Let \(I\) be an internally injective \(R\)-module in \(\Xc\), and let \(X \in \Xc\).
  If \(X\) is covered by a set, then \(X \times I\) is an internally injective \((X{\times}R)\)-module.
\end{pro}

\begin{proof}
  Let \(I\) be an internally injective \(R\)-module in \( \Xc \), and let \( X \in \Xc \).
  We wish to show that \( X \times I \) is an internally injective \((X{\times}R)\)-module in \( \Xc / X \).
  Since \(X\) is covered by a set, we have an effective epimorphism \(V \thra X\) with \(0\)-truncated domain.
  By \cref{lem:internal-injectivity-stable-sets}, \(V \times I\) is an internally injective \((V{\times}R)\)-module.
  But \(V \thra X \) is a \((-1)\)-connected object over \(X\), thus by \cref{lem:internal-injectivity-descends} we can descend internal injectivity from \(V \times I \) to \(X \times I\).
\end{proof}

\begin{cor} \label{cor:sets-cover-internal-injectivity-stable}
  If sets cover in \( \Xc \), then internal injectivity of \(R\)-modules is stable by base change.
  \qed
\end{cor}

Our next goal is to extend this result to any slice of an \(\infty\)-topos in which sets cover.
This generalization will let us understand the interpretation of internal injectivity when working in a non-empty context in HoTT, which corresponds to working in a slice of the chosen \(\infty\)-topos model.
The key result is the following:

\begin{pro} \label{pro:n-types-cover-slice}
  If \(n\)-types cover in \(\Xc\), then \(n\)-types cover in \(\Xc / X\) for any \((n+1)\)-type \(X \in \Xc\).
\end{pro}

\begin{proof}
  Let $X$ be in $(n+1)$-type in $\Xc$ and
  consider an object \(Y \to X\) in \(\Xc/X\).
  Since \(n\)-types cover in \(\Xc\), there is an effective epimorphism \( e : V \thra Y\) with \(V\) an \(n\)-type in \(\Xc\).
  The map \(e\) defines an effective epimorphism over \(X\) with domain the composite \(V \to Y \to X\).
  The latter is a map from an \(n\)-type to an \((n+1)\)-type in $\Xc$, which is necessarily \(n\)-truncated (as is easily shown in HoTT, by showing that the fibers are all $n$-types).
  Hence the domain of \(e\) is \(n\)-truncated as an object of \(\Xc / X\), so \(Y\) is covered by an \(n\)-type.
\end{proof}

\begin{thm} \label{thm:sets-cover-stability}
  Suppose sets cover in \( \Xc \), and let \( X \in \Xc \).
  For any ring \( R \in \Xc \kslice X \), internal injectivity of \( R \)-modules is stable by base change in \( \Xc / X \).
\end{thm}

\begin{proof}
  Let \( I \in \Xc / X \) be an internally injective \(R\)-module.
  The truncation map \( X \to \Tr{X}_1 \) is $1$-connected, and therefore induces an equivalence \( \sets{\Xc / \Tr{X}_1} \xra{\sim} \sets{\Xc / X} \) by base change~\cite[Lemma~7.2.1.13]{HTT}.
  Moving back along this equivalence, $I$ becomes an internally injective module over \( \Tr{X}_1 \).
  Since \(\Tr{X}_1\) is a 1-type, the previous proposition implies that sets cover in \(\Xc / \Tr{X}_1\).
  By \cref{cor:sets-cover-internal-injectivity-stable}, this means \(I\) is internally injective when pulled back to any slice of \( \Xc / \Tr{X}_1 \).
  But any slice of \( \Xc / X \) is a slice of \( \Xc / \Tr{X}_1 \), so we are done.
\end{proof}

\begin{rmk}
  These methods apply in much greater generality than just internal injectivity.
  Consider, for example, some ``internal property'' \(P\) of an object (or structure) in a \(1\)-topos which is stable by base change and descends along (effective) epimorphisms.
  Then, since the $0$-truncated objects in $\Xc$ form a $1$-topos,
  we can ask whether the property \(P\) holds for some given \(0\)-truncated object \( Y \) in any slice \( \Xc / X\).
  The arguments above show that if sets cover in \( \Xc \),
  then the property \( P \) is stable by base change in \( \Xc / X\).
 (Making precise the meaning of ``internal property'' is beyond our current scope.)
\end{rmk}

\begin{cor} \label{cor:hott-internal-injectivity}
  Suppose sets cover in \( \Xc \).
  Consider an object \( X \), a ring \( R \in \Xc \kslice X \)
  and an \( R \)-module \( I \).
  Then $I$ is HoTT-injective if and only if it is internally injective. \qed
\end{cor}

Using \cref{cor:hott-internal-injectivity}, we explain how the interpretation of our Ext groups from HoTT recover the classical notion of sheaf Ext.

\begin{dfn}\label{dfn:sheaf-Ext}
Let $\Ec$ be a $1$-topos, let $R$ be a ring object in $\Ec$,
and let $B$ be an $R$-module.
We define the functor $\sExt{R}^n(B, -) : \eMod{R} \to \eAb{\Xc}$ to be the
$n^{th}$ right derived functor of $\ihMod{R}(B, -)$,
where we use an external injective resolution to define the derived functor.
We refer to $\sExt{R}^n(B, A)$ as \define{sheaf Ext}.
We extend this definition to an $\infty$-topos $\Xc$ by applying it to the
\(1\)-topos \( \sets{\Xc} \).
\end{dfn}

The sheaf Ext groups arise in algebraic geometry~(\cite[Chapitre~IV]{Gro57},
\cite[Section~III.6]{hartshorne})
and are also considered in~\cite[Section~18.4]{KS06} and~\cite[Section~13.4]{Ble17}.

For any \(R\)-module \(B\) in \(\Xck\), we obtain an internal functor \(\iExt{R}^n(B,-)\) in \(\Xc\) by interpretation, which yields an ordinary functor
\( \iExt{R}^n(B,-) : \eMod{R} \to \eAb{\Xc}' \).
Here the \('\) indicates that the codomain consists of abelian group objects in \(\Xc\), not just \(\Xck\).
The dual of \cref{pro:ext-projective-resolution} lets us compute \( \iExt{R}^n(B,A) \) via a HoTT-injective resolution of \( A \).
Combining the results of this section, we obtain the following:

\begin{thm} \label{thm:recover-sheaf-ext}
  Suppose sets cover in \( \Xc \).
  For any \( X \in \Xc \), ring \( R \in \Xc \kslice X \) and \( R \)-module \( B \), the functor \( \iExt{R}^n(B,-) : \eMod{R} \to \eAb{\Xc / X}' \) is naturally isomorphic to the sheaf Ext functor $\sExt{R}^n(B,-)$.
  In particular, we may take \(\iExt{R}^n(B,-)\) to land in \(\eAb{\Xc / X}\).
\end{thm}

\begin{proof}
  Since (external) injectives in \(\eMod{R}\) are always internally injective by~\cite[Theorem~3.8]{Ble18}, and moreover internal and HoTT-injectives coincide in \( \Xc / X \) by \cref{cor:hott-internal-injectivity}, we can use an (externally) injective resolution to compute \(\iExt{R}^n(B,-)\) by the dual of \cref{pro:ext-projective-resolution}.
  But this means that \(\iExt{R}^n(B,-)\) is the \(n^{th}\) right derived functor of the internal hom of \(R\)-modules, meaning it is naturally isomorphic to \(\sExt{R}^n(B,-)\).
  In particular, it is \(\kappa\)-compact.
\end{proof}

It follows that the computation in \cref{pro:ext-computation} can be regarded
as a computation of sheaf Ext in the Sierpi\'nski \(\infty\)-topos,
using a HoTT-projective resolution.

\begin{rmk}
In his thesis, Blechschmidt gives a definition of sheaf Ext groups in the internal language of a localic \(1\)-topos using the existence of enough injectives~\cite[Section~13.4]{Ble17}.
In contrast, our internal Ext is the interpretation of \cref{dfn:ext-1}, which does not rely on injectives.
\end{rmk}

Since there are always enough internally injective \(R\)-modules (\cref{rmk:enough-internal-injectives}), we deduce the following:

\begin{cor}\label{cor:global-points-sheaf-Ext}
  Let \(X \in \Xc\).
  Suppose that sets cover in \(\Xc\) and that the (set-restricted) global points functor \(\Gamma_X : \eSet{\Xc/X} \to \eSet{} \) preserves epimorphisms.
  For any ring \(R \in \Xc \kslice X\), \(R\)-modules \(B\) and \(A\), and \(n \geq 0\), we have a natural isomorphism
  \( \Gamma_X \iExt{R}^n(B,A) \ \simeq \ \eExt{R}^n(B,A) .\)
  In particular, the ordinary Ext groups are obtained as the global points of sheaf Ext.
\end{cor}

\begin{proof}
  An argument similar to \cref{pro:internal-implies-external-projectivity} shows that internal injectivity implies external injectivity of \(R\)-modules under our assumption on \(\Gamma_X\).
  Thus internal and external injectivity coincide, and are equivalent to HoTT-injectivity by \cref{cor:hott-internal-injectivity}.
  The statement follows by the same proof as in \cref{pro:global-points-HoTT-ext-projective-resolution}, but using an (internally) injective resolution of \(A\).
\end{proof}

\begin{rmk}
  It is well-known that the global points of sheaf Ext recover the ordinary Ext groups whenever the global points functor \(\Gamma_X\) preserves epimorphisms (hence is exact).
  This fact is an easy consequence of the (``local-to-global'') Grothendieck spectral sequence which relates sheaf Ext and ordinary Ext, specifically:
  \[ (\on{R}^p \Gamma) \sExt{R}^q(B,A) \implies \eExt{R}^{p + q}(B,A) . \]
  (Here \(\on{R}^p\) denotes the \(p^{th}\) right derived functor.)
  Our assumption on \(\Gamma_X\) implies that \(\on{R}^p \Gamma_X\) vanishes for \(p > 0\), which means this spectral sequence collapses at the \(E_2\)-page.
  It immediately follows that we have an isomorphism \(\Gamma_X \sExt{R}^n(B,A) \simeq \eExt{R}^n(B,A)\), for all \(n \in \Nb\).
\end{rmk}

We also record the following corollary of \cref{thm:recover-sheaf-ext}.

\begin{cor} \label{cor:ext-determined-by-sets}
  Suppose sets cover in \( \Xc \), and let \( X \in \Xc \).
  The interpretation of \( \Ext^n_R(B,A) \) into \( \Xc / X \) depends solely on \( \sets{\Xc / X} \). \qed
\end{cor}

There are many \( \infty \)-toposes which share $1$-toposes of $0$-truncated objects.
For example, if \( X \in \cat{S} \) is a pointed, connected space, then $0$-truncated objects in the slice \( \infty \)-topos \( \cat{S} / X \) are \( \pi_1(X) \)-sets.
Thus if \( X \) is simply connected, these are just sets.
Since sets cover in \( \spaces \), the corollary tells us that interpreting \( \iExt{R}^n \) into any slice \( \spaces / X \) with \( X \) simply connected yields the same result (up to equivalence).
This means in particular that we can move between \( \spaces \) and \( \spaces / X \) when computing \( \iExt{R}^n \)---a potentially useful trick.

\subsection{Ext over \texorpdfstring{\(BG\)}{BG}} \label{ssec:Ext-BG}

Let \(X\) be a pointed, connected object in \(\Xc\), with base point \(\pt : X\).
In this final section, we study Ext groups of abelian group objects in the slice $\Xc / X$, and relate them to Ext groups in the base \(\Xc\).
As we will see, these considerations are intimately related with those of \cref{ssec:Ext-ZG}, and they illustrate the theory developed thus far in \cref{sec:interpretation}.

We refer to abelian group objects in $\Xc / X$ as $(X {\times} \Zb)$-modules,
as this makes the base clear.
We mention that the 1-truncation map \(X \to B \pi_1(X)\) is 1-connected and therefore induces an equivalence between the 1-topos of sets over \(X\) and the 1-topos of sets over \(B \pi_1(X)\) by pulling back~\cite[Lemma~7.2.1.13]{HTT}.
Accordingly, we get an equivalence \(\eMod{X \times \Zb} \simeq \eMod{B \pi_1(X) \times \Zb}\) of categories of modules.
To emphasize that no truncation assumptions are needed, we work with \(X\) rather than \(B \pi_1(X)\).

The category of \((X{\times}\Zb)\)-modules has another description.
However, as we will see in a moment, this other description is not equivalent when working internally, since it changes the ambient topos.
For any (0-truncated) group object $G$ in \(\Xc\), we can form the internal group ring \(\Zb G\) as the object \(\bigoplus_G \Zb\) with its natural ring structure.
This is the result of interpreting the group ring of \cref{cst:zg} into \(\Xc\).
Being a \(G\)-set, \(\Zb G\) defines an object of \(\Xc / BG\) which may be seen to be the free abelian group on the base point (or universal cover) \(1 \to BG\), though we will not make use of this description.

\begin{pro} \label{pro:ZG-mod}
  Restriction to the base point of $X$ gives an equivalence of \(1\)-categories
  \[ \eMod{X \times \Zb} \ \simeq \ \eMod{\Zb \pi_1(X)} , \]
  where the left-hand side is the category of abelian groups in \(\Xc \kslice X\) and the right-hand side is the category of \(\Zb \pi_1(X)\)-modules in~\(\Xck\).
  In particular, we obtain an isomorphism
  \[ \eExt{X \times \Zb}^n(B, A) \ \simeq \ \eExt{\Zb \pi_1(X)}^n(B_{\pt}, A_{\pt}) \]
 of external Ext groups for all $(X {\times} \Zb)$-modules $A$ and $B$.
\end{pro}

We point out that this statement does \emph{not} imply that \(\iExt{X \times \Zb}\) coincides with \(\iExt{\Zb \pi_1(X)}\), as the former is an abelian group in \(\Xc / X\), whereas the latter is an abelian group in \(\Xc\).
However, the relation between these objects is interesting and is further discussed below.

\begin{proof}
  By~\cite[Theorem~4.3.4]{Fla22}, the category \(\eMod{X \times \Zb}\) is equivalent to the category \(\Xc(X, \iMod{\Zb})\) (whose categorical structure comes from \(\iMod{\Zb}\)).
  The interpretation of \cref{pro:zg-mod} yields an equivalence of categories
  \( \uXc(X, \iMod{\Zb}) \simeq \iMod{\Zb \pi_1(X)} \) in \(\Xc\),
  where $\uXc(X, \iMod{\Zb})$ denotes the internal category whose object of objects
  is the internal hom in $\Xc$.
  On global points this yields the desired equivalence of categories
  \( \Xc(X, \iMod{\Zb}) \simeq \eMod{\Zb \pi_1(X)} \).
  It follows that the stated (external) Ext groups are isomorphic.
\end{proof}

Given a $(X {\times} \Zb)$-module $A$, we call its restriction $A_{\pt}$ along $1 \to X$
the \define{underlying abelian group object} of $A$.
This has a natural $\pi_1(X)$-action, so it can also be regarded as a $\Zb \pi_1(X)$-module in $\Xc$.
Note that the equivalence \( \eMod{X \times \Zb} \simeq \eMod{\Zb \pi_1(X)} \)
sends the ring $X \times \Zb$ to $\Zb$ with the trivial $\pi_1(X)$-action,
and not to the ring $\Zb \pi_1(X)$.
We also warn the reader that care must be taken when moving across this equivalence.
For example, the category \(\Mod{X \times \Zb}\) is enriched over itself via the internal hom of abelian groups in \(\Xc \kslice X\),
while \(\Mod{\Zb \pi_1(X)}\) is naturally enriched over \(\eAb{\Xc}\).
These hom-objects live in different categories, but one can check that the latter is the \(\pi_1(X)\)-fixed points of the former:
\[ \prod_{X} \ihMod{X{\times}\Zb}(B,A) \ \simeq \ \ihMod{\Zb \pi_1(X)}(B_*, A_*) . \]
Because of the difference between these internal homs, many internal properties are not preserved by the equivalence of \cref{pro:ZG-mod}.
An example of this is given in \cref{exa:internal-not-external-projective}, as explained in the discussion immediately after it.

We record the following fact, which immediately follows from base-change stability of interpretation.

\begin{pro} \label{pro:zg-ext-slice}
  Let \(B\) and \(A\) be \((X{\times }\Zb)\)-modules.
  The underlying abelian group object of the \((X{\times}\Zb)\)-module \(\iExt{X \times \Zb}^n(B,A)\) is \(\iExt{\Zb}^n(B_{\pt},A_{\pt})\). \qed
\end{pro}

A concrete description of the \(\pi_1(X)\)-action on \(\iExt{\Zb}^n(B_\pt, A_\pt)\) can be worked out from \cref{pro:zg-ext-description} and the discussion surrounding it.

\begin{rmk}\label{rmk:fixed-points}
  It might be tempting to believe that the abelian group object \(\iExt{\Zb \pi_1(X)}^n(B_{\pt},A_{\pt})\) is isomorphic to the fixed points of the \(\Zb \pi_1(X)\)-module \(\iExt{X \times \Zb}^n(B,A)_{\pt}\) described by the previous theorem.
  In general, this is not the case, as we will see in \cref{exa:ext-G-modules}.
\end{rmk}

We deduce a vanishing result for \(\iExt{X \times \Zb}\).

\begin{cor}\label{cor:zg-ext-vanish}
  Let \(n\) be a natural number.
  Suppose that \(\iExt{\Zb}^n(B,A)\) vanishes for all abelian groups \(B\) and \(A\) in \(\Xck\).
  Then \(\iExt{X \times \Zb}^n(N,M)\) also vanishes for all \((X {\times} \Zb)\)-modules \(N\) and \(M\). \qed
\end{cor}

Our next result characterizes internal injectivity of abelian groups in the slice \(\Xc / X\), and generalizes~\cite[Proposition~1.2(i)]{Har81} for ordinary sheaves on a space.

\begin{pro} \label{pro:internally-injective-G-module}
  An \((X{\times}\Zb)\)-module is internally injective if and only if its underlying abelian group object is internally injective.
  The same holds for internal projectivity.
\end{pro}

\begin{proof}
  We prove the injective case, as the projective case is shown similarly.

  \((\rightarrow)\) Let \(I\) be an internally injective \((X{\times}\Zb)\)-module.
  A monomorphism  \(i : A \hra B\) of abelian groups in \(\Xc\) pulls back to a monomorphism \(i : X \times A \hra X \times B\) between \((X{\times}\Zb)\)-modules (with trivial action).
  Thus we get an epimorphism \(i^* : \ihMod{X \times \Zb}(X \times B, I) \thra \ihMod{X \times \Zb}(X \times A, I)\) of \((X{\times}\Zb)\)-modules.
  This homomorphism is given by \(i^* : \ihMod{\Zb}(B,I_\pt) \to \ihMod{\Zb}(A,I_\pt) \) on the underlying abelian group objects, since base change (here along \(1 \thra X\)) respects internal homs.
  The latter map is therefore an epimorphism, as desired.

  \((\leftarrow)\) By \cref{lem:internal-injectivity-descends} we can descend internal injectivity along the effective epimorphism \(1 \thra X\), meaning \(I\) is an internally injective \((X{\times}\Zb)\)-module whenever \(I_*\) is an internally injective abelian group object.
\end{proof}

We note that the proof of \((\rightarrow)\) only used that \(X\) was pointed.

\begin{rmk}
The previous proposition gives another way of understanding \cref{pro:zg-ext-slice}.
Namely, if one computes \(\iExt{X{\times}\Zb}^n(B,A)\) using an internally injective resolution of \((X{\times}\Zb)\)-modules (which always exists), then the underlying abelian resolution can be used to compute \(\iExt{\Zb}^n(B_\pt, A_\pt)\).
Thus we see that the latter is the underlying abelian group object of the former.
\end{rmk}

Before our next examples, we show that internal projectivity and HoTT-projectivity coincide in spaces.

\begin{pro}\label{pro:hott-projectivity-spaces}
  In the \(\infty\)-topos of spaces, HoTT-projectivity and internal projectivity of modules coincide.
\end{pro}

\begin{proof}
  Firstly, note that external and internal projectivity coincide in spaces.
  Now, suppose that \(P\) is an (internally) projective \(R\)-module.
  By \cref{pro:hott-projectivity}, we need to show that \(X \times P\) is an internally projective \((X{\times}R)\)-module in \(\Xc / X\), for any \(X \in \Xc\).
  Since sets cover in spaces, and internal projectivity descends along effective epimorphisms by a proof analogous to \cref{lem:internal-injectivity-descends}, we can assume that \(X\) is a set.
  Then an \((X{\times}R)\)-module is simply an \(X\)-indexed collection of \(R\)-modules, and the internal hom of such is the indexwise hom.
  The axiom of choice implies that an \(X\)-indexed collection of epimorphisms defines an epimorphism between the collections, so \(X \times P\) is internally projective.
\end{proof}

We now give examples of modules which are internally projective, but not externally projective.

\begin{exa} \label{exa:internal-not-external-projective}
  Let \(G\) be a non-trivial (0-truncated) group in \(\spaces\).
  The \((BG{\times}\Zb)\)-module \(BG \times \Zb\) is internally projective, but not externally projective.
  It follows that there exists a $(BG{\times}\Zb)$-module $A$ so that
  $\iExt{BG{\times}\Zb}^1(BG \times \Zb, A) = 0$ and $\eExt{BG{\times}\Zb}^1(BG \times \Zb, A) \neq 0$.
\end{exa}

\begin{proof}
  Any ring is HoTT-projective as a module over itself, and is therefore internally projective.
  In particular, \(BG \times \Zb\) is internally projective.
  (This can also be seen from \cref{pro:internally-injective-G-module}.)
  External projectivity of \((BG{\times}\Zb)\)-modules corresponds to ordinary projectivity of \(\Zb G\)-modules by \cref{pro:ZG-mod}.
  As a \(\Zb G\)-module, \(BG \times \Zb\) corresponds to the abelian group \(\Zb\) with trivial \(G\)-action.
  Since \(G\) is non-trivial, the augmentation homomorphism \(\Zb G \to \Zb\) cannot split, thus \(\Zb\) is not projective.
  The last claim follows from \cref{pro:HoTT-proj}.
\end{proof}

As mentioned above, for a module in $\spaces$, internal and external projectivity agree.
So the example shows that while $BG \times \Zb$ is internally projective as a $(BG {\times} \Zb)$-module,
the corresponding $\Zb G$-module $\Zb$ is not internally projective.
This demonstrates the sense in which the equivalence of \cref{pro:ZG-mod} does not respect internal properties, since the ambient topos changes.
This is also demonstrated by the following example, which ties together many of our results and remarks into a concrete example in spaces.

\begin{exa} \label{exa:ext-G-modules}
  Take $\Xc$ to be $\spaces$ and \(G\) to be a (0-truncated) group.
  Since \(\eAb{\spaces} = \eAb{}\) has global dimension 1 and $\iExt{\Zb}^n$ interprets to ordinary $\eExt{\Zb}^n$ by \cref{pro:recover-Ext-spaces}, we deduce that $\iExt{\Zb}^n$ vanishes for $n > 1$.
  \cref{cor:zg-ext-vanish} then says that \(\iExt{BG \times \Zb}^n\) vanishes for $n > 1$.
  Even more, $\iExt{BG \times \Zb}^n(BG \times \Zb, M)$ vanishes for all $n \geq 1$
  and every $M$, since $BG \times \Zb$ is internally projective.
  On the other hand, by \cref{pro:ZG-mod},
  the ordinary Ext groups $\eExt{BG \times \Zb}^n(BG \times \Zb, A)$ are the same as
  the ordinary Ext groups $\eExt{\Zb G}^n(\Zb, A_{\pt})$, which need not vanish.
  For example, it is well known that \( \eExt{\Zb G}^n (\Zb, M) \cong H^n(BG; M) \), which may be nonzero for all \(n \in \Nb \).
  Indeed, as we saw in \cref{exa:internal-not-external-projective}, \(\Zb\) with trivial action is not a projective \(\Zb G\)-module.
  Note also that $\iExt{\Zb G}^n$ agrees with $\eExt{\Zb G}^n$,
  again using \cref{pro:recover-Ext-spaces}.
  In particular, as mentioned in \cref{rmk:fixed-points}, it is not the case in
  general that $\iExt{\Zb G}^n$ can be described as the fixed points of the \(\Zb G\)-module corresponding to \(\iExt{BG \times \Zb}^n\), even for $n = 1$.

  We explain this phenomenon in a bit more detail.
  A short exact sequence of \((BG{\times}\Zb)\)-modules
  \[ 0 \to A \to E \to B \to 0 \]
  may be seen both as an element of the abelian group \(\eExt{BG \times \Zb}^1(B, A) \cong \eExt{\Zb G}^1(B_\pt,A_\pt)\), and as an element (indeed, \(G\)-fixed point) of the \(\Zb G\)-module \(\iExt{BG \times \Zb}^1(B,A)_{\pt}\).
  This extension \(E\) is trivial as an element of the former if and only if the epimorphism \(E \to B\) admits a \((BG{\times}\Zb)\)-module section (i.e., a \(G\)-equivariant section).
  In contrast, \(E\) is trivial as an element of the latter if and only if the underlying homomorphism \(E_\pt \to B_\pt\) of abelian groups admits an \emph{abelian} section, by \cref{pro:zg-ext-slice}.
\end{exa}

\printbibliography

\end{document}